\documentclass[12pt]{article}

\usepackage{amsthm,amsmath,amsfonts,epsfig,amssymb,latexsym}
\usepackage{mathrsfs}  
\usepackage[T1]{fontenc}
\usepackage{graphicx}
\usepackage{enumerate}
\usepackage{xy,xypic}
\usepackage{graphics}
\usepackage{footnote}
\usepackage{scalefnt}
\usepackage{tikz}
\usetikzlibrary{shapes,arrows,fit,calc,positioning}
\usetikzlibrary{matrix}
\xyoption{all}
\title{{\Large Derived Witt-D\'{e}vissage Formalism}} 
\author{
 Satya Mandal \footnote{Partially supported by a General Research Grant from 
 KU}\\ 
{\small  University of Kansas, Lawrence KS 66045} \\
{\small {\it email: mandal@ku.edu}} \\
}

\begin{document}
\renewcommand{\baselinestretch}{1.255}
\setlength{\parskip}{1ex plus0.5ex}
\date{January 27, 2014}
\newtheorem{theorem}{Theorem}[section]
\newtheorem{proposition}[theorem]{Proposition}
\newtheorem{lemma}[theorem]{Lemma}
\newtheorem{corollary}[theorem]{Corollary}
\newtheorem{construction}[theorem]{Construction}
\newtheorem{notations}[theorem]{Notations}
\newtheorem{question}[theorem]{Question}
\newtheorem{example}[theorem]{Example}
\newtheorem{definition}[theorem]{Definition}
\newtheorem{remark}[theorem]{Remark}

\newcommand{\iso}{\stackrel{\sim}{\longrightarrow}}
\newcommand{\sur}{\twoheadrightarrow}
\newcommand{\bD}{\begin{definition}}
\newcommand{\eD}{\end{definition}}
\newcommand{\bP}{\begin{proposition}}
\newcommand{\eP}{\end{proposition}}
\newcommand{\bL}{\begin{lemma}}
\newcommand{\eL}{\end{lemma}}
\newcommand{\bT}{\begin{theorem}}
\newcommand{\eT}{\end{theorem}}
\newcommand{\bC}{\begin{corollary}}
\newcommand{\eC}{\end{corollary}}
\newcommand{\eop}{\hfill \square\\}
\newcommand{\pf}{\noindent{\bf Proof.~}}
\newcommand{\PD}{\text{proj} \dim}
\newcommand{\lra}{\longrightarrow}
\newcommand{\hra}{\hookrightarrow}
\newcommand{\Lra}{\Longrightarrow}
\newcommand{\Llra}{\Longleftrightarrow}
\newcommand{\bE}{\begin{enumerate}}
\newcommand{\eE}{\end{enumerate}}
\newcommand{\pic}{The proof is complete.}
\newcommand{\TCP}{\textcolor{purple}}
\newcommand{\TCM}{\textcolor{magenta}}
\newcommand{\TCR}{\textcolor{red}}
\newcommand{\TCB}{\textcolor{blue}}
\newcommand{\TCG}{\textcolor{green}}

\def\spec#1{\mathrm{Spec}(#1)}
\def\m{\mathfrak {m}}
\def\CA{\mathcal {A}}
\def\CB{\mathcal {B}}
\def\CP{\mathcal {P}}
\def\CC{\mathcal {C}}
\def\CD{\mathcal {D}}
\def\CE{\mathcal {E}}
\def\CF{\mathcal {F}}
\def\CG{\mathcal {G}}
\def\CH{\mathcal {H}}
\def\CI{\mathcal {I}}
\def\CJ{\mathcal {J}}
\def\CK{\mathcal {K}}
\def\CL{\mathcal {L}}
\def\CM{\mathcal {M}}
\def\CN{\mathcal {N}}
\def\CO{\mathcal {O}}
\def\CP{\mathcal {P}}
\def\CQ{\mathcal {Q}}
\def\CR{\mathcal {R}}
\def\CS{\mathcal {S}}
\def\CT{\mathcal {T}}
\def\CU{\mathcal {U}}
\def\CV{\mathcal {V}}
\def\CW{\mathcal {W}}
\def\CX{\mathcal {X}}
\def\CY{\mathcal {Y}}
\def\CZ{\mathcal {Z}}

\def\SI{\mathscr {I}}
\def\SV{\mathscr {V}} 

\newcommand{\smallcirc}[1]{\scalebox{#1}{$\circ$}}
\def\BA{\mathbb {A}}
\def\BB{\mathbb {B}}
\def\BC{\mathbb {C}}
\def\BD{\mathbb {D}}
\def\BE{\mathbb {E}}
\def\BF{\mathbb {F}}
\def\BG{\mathbb {G}}
\def\BH{\mathbb {H}}
\def\BI{\mathbb {I}}
\def\BJ{\mathbb {J}}
\def\BK{\mathbb {K}}
\def\BL{\mathbb {L}}
\def\BM{\mathbb {M}}
\def\BN{\mathbb {N}}
\def\BO{\mathbb {O}}
\def\BP{\mathbb {P}}
\def\BQ{\mathbb {Q}}
\def\BR{\mathbb {R}}
\def\BS{\mathbb {S}}
\def\BT{\mathbb {T}}
\def\BU{\mathbb {U}}
\def\BV{\mathbb {V}}
\def\BW{\mathbb {W}}
\def\BX{\mathbb {X}}
\def\BY{\mathbb {Y}}
\def\BZ{\mathbb {Z}}

\def\LF{\mathcal {LF}}
\maketitle
\noindent{\bf Abstract}: {\it 
In this article we establish some formalism of  Derived 
Witt-D\'{e}vissage theory for resolving subcategories of abelian categories. 
Results directly apply to noetherian schemes. 
}
\section{Introduction}
This is another one of a series of articles
(\cite{MSd, MSfpd, MSf}) dedicated to study  of derived Witt groups.
Our interest in derived Witt groups  emanates from 
our interest in Chow-Witt groups,
as  obstruction groups for projective modules.
We would refrain from repetition of the 
 introductory comments and the 
 background discussion provided in (\cite{MSd}).

In (\cite{MSd}), we proved  d\'{e}vissave theorems 
for Cohen-Macaulay affine schemes,  extending the
corresponding theorem of Balmer and Walter (\cite{BW}) 
on regular schemes.  
In this paper, we prove analogous d\'{e}vissave theorems
for  noetherian schemes $X$, with the hypothesis that 
any coherent sheaf on $X$ is a quotient of a locally free
sheaf on $X$. In particular, we remove the Cohen-Macaulay
condition in (\cite{MSd}).
In the process of doing so, it became clear that
the arguments could be made  formal enough,
so that they work for resolving subcategories 
of abelian categories. The concept of resolving subcategories
dates back to the paper of Auslander and  Bridger (\cite{AB}).
More recently, there has been considerable amount of activities  
({\it e.g.} \cite{T1,T2})  
on resolving subcategories of 
the category $Mod(A)$ of finitely generated modules over
noetherian commutative rings $A$. 
Much of it is directed toward the classification of such resolving 
subcategories under various conditions, which is encompassed by 
similar attempts to classify variety of types 
of subcategories of the module categories. 
The work of Benson, Iyengar and Krause  ({\it e.g.} \cite{BIK})
would be one of the stimulus, 
where they consider subcategories of group algebras. 
To the best of our knowledge,  very little literature is
available on resolving subcategories of abelian categories.
We give a version of the d\'{e}vissage theorem for 
resolving subcategories
of abelian categories, 
which directly applies to noetherian schemes.  
From our perspective, 
the results on Witt groups of noetherian schemes constitute  
the essence of this article.
However, the formal version of the results 
on resolving subcategories unifies the theory and
has its significance by its own rights. 
In this introduction, we will 
state the versions of the results for resolving subcategories,
as follows.

Suppose ${\SV}$ is a resolving subcategory (see definition \ref{resCatDef})
of an abelian category ${\CC}$. Let $\omega$ be an object
in ${\CC}$ with an injective resolution. For objects 
$M$ in ${\CC}$, denote $M^*:=Hom(M, \omega)$.  Assume 
 ${\SV}$ inherits a duality structure from $^*$ and 
 ${\SV}$ is totally $\omega$-reflexive (see definition \ref{defResWduality}). 
The following is a list of a few notations and facts. 
\bE
\item Let ${\CB}({\SV}, \omega)$ denote the full subcategory of 
objects in ${\CC}$ with finite ${\SV}$-dimension. 
Assume $d:=\max \{\dim_{{\SV}}M: M\in {\CB}({\SV}, \omega) \}<\infty$.
\item Let 
$
{\CA}(\omega):=\{M\in {\CB}({\SV}, \omega): 
Ext^i(M,\omega)=0~~\forall~i<d\}.
$ 
It follows that $M\mapsto M^{\vee}:=Ext^d(M, \omega)$ is a 
duality on ${\CA}(\omega)$.
\item 
$D^b({\SV})$ would denote the bounded derived category of 
complexes of objects in ${\SV}$ and $D^b_{{\CA}(\omega)}({\SV})$
would be the subcategory of complexes with homologies in 
${\CA}(\omega)$.
\item It is a standard fact (see \ref{formalZeta})
that there is a canonical functor
$\zeta: {\CB}({\SV}, \omega) \lra D^b({\SV})$,
given by (choices of)  ${\SV}$-resolutions.
\item The restriction of $\zeta:{\CA}(\omega)\lra D^b_{{\CA}(\omega)}({\SV})$
is a duality preserving  functor. 
\eE 
First, we prove isomorphisms of Witt groups, as follows. 
$$
W({\CA}(\omega), ^{\vee}, \pm \tilde{\varpi})\iso
W(D^b_{{\CA}(\omega)}({\CA}(\omega), ^{\vee}, \pm \tilde{\varpi}))\iso
 W(D^b({\CA}(\omega), ^{\vee}, \pm \tilde{\varpi})).\\ 
 $$
It follows from a theorem of Balmer (\cite{TWGII}) that the 
first and the last groups are isomorphic.
We further  prove that the functor $\zeta$ induces isomorphims
of Witt groups
$$
W_{St}^+({\CA}(\omega)) \iso 
W^d\left(D^b_{{\mathcal A}(\omega)}({\SV}), *, 1, \varpi \right),
\quad 
W_{St}^-({\CA}(\omega)) \iso W^{d-2}\left(D^b_{{\mathcal A}(\omega)}
({\SV}), *, 1, -\varpi \right)   
$$
where subscript "St" corresponds to "standard"
sign convention of the duality.
Also for $n=d-1, d-3$, we have
$W^{n}\left(D^b_{{\mathcal A}(\omega)}({\SV}, *,1, \pm \varpi)\right)=0$.
By 4-periodicity all the shifted Witt groups are determined. 

These results can be applied to 
noetherian schemes $X$, with ${\CC}=Coh(X)$ and ${\SV}$ as the 
subcategory of locally free sheaves on $X$, provided ${\SV}$
is a resolving subcategory. This will be the case for a wide variety 
of schemes $X$ (see \cite{H}), 
including those that have an ample invertible sheaf. In these applications,
we assume 
$d:=\dim X=\max\{depth({\CO}_{X,x}): x~{\rm is~a~closed~point} \}$ 
(see remark \ref{lowDepth}). 
As a consequence, the following decomposition theorem follows.
\bT 
Suppose $X$ is a noetherian scheme, with $\dim X=d$, as in
(\ref{notaCate}) and $X^{(d)}$ will denote the
set of all closed points of codimension $d$ in $X$. We assume 
$d=\max\{depth({\CO}_{X,x}): x~{\rm is~a~closed~point} \}=\dim X$.  
Then, the  homomorphisms
$$
W^d\left(D^b_{{\mathcal A}(X)}({\SV}(X)), *, 1, \varpi \right)\iso
\bigoplus_{x\in X^{(d)}} 
W^d(D^b_{{\CA}({\CO}_{X,x})}({\CO}_{X,x}),*,1,\varpi) 
\qquad {\rm and} 
$$
$$
W^{d-2}\left(D^b_{{\mathcal A}(X)}({\SV}(X)), *, 1, -\varpi \right)\iso
\bigoplus_{x\in X^{(d)}} 
W^{d-2}(D^b_{{\CA}({\CO}_{X,x})}({\CO}_{X,x}),*,1,-\varpi) 
$$
are isomorphims.
\eT
When $X$ is regular,  this is a theorem of Balmer and Walter (\cite{BW}).
In this case,
first  isomorphism
 follows from such a decomposition of the corresponding
 categories
 and the second isomorphism would have zeros on both sides.    

Before we conclude this introduction, we comment on the sense of direction of this series of articles (\cite{MSd, MSfpd, MSf}). 
For  noetherian schemes $X$ as above, 
we will give a Gersten-Witt like complex of the "relative" Witt groups 
$W^p({\CD}^p(X), {\CD}^{p+1}(X))$, where ${\CD}^p(X)$ denotes the subcategory of  the finite derived category 
$D^b(X)$, of complexes 
${\CE}_{\bullet}$ with  finite locally free dimension homologies 
${\CH}_i({\CE}_{\bullet})$ and 
$co\dim({\CH}_i({\CE}_{\bullet}))\geq p$. The "relative" Witt groups are defined by forming 
a group of isometry classes of $S$-spaces, as in (\cite{TWGI}). These results will appear subsequently.

We close this introduction,
with a few comments on the lay out of this article.
As would be expected,
 the results on Witt groups of noetherian schemes would follow
 from that of the resolving subcategories.
 However, we made a choice to
 give complete proofs for the former (\S~\ref{mainIsoSEC},~\ref{finalResSec})
 and the proofs of the  results on resolving subcategories
(\S ~\ref{FormalDiv}) would follow similarly.
Throughout in \S~\ref{mainIsoSEC},~\ref{finalResSec}, we remind the readers
that our proofs are formal enough to apply to \S ~\ref{FormalDiv}.
We also give a proof of the existence of the functor 
$\zeta: {\CB}(X) \lra D^b({\SV}(X))$ in \S ~\ref{zetaSection}. While
this is a standard fact among the experts, we were unable to find
 references that would be of any help for the readership of this article. 
While the main isomorphism theorem was dealt with in \S ~\ref{mainIsoSEC}
, we summarize 
our results on noetherian schemes in \S ~\ref{finalResSec}.
In \S ~\ref{FormalDiv}, we state our results on
resolving subcategories. 

\vspace{2mm}
\noindent{\bf Acknowledgements}: 
{\it I would like to sincerely thank Sarang Sane for numerous discussions I had with him for a long period of time. I also thank
 him for  suggesting improvements of some 
 of the proofs. 
 Thanks are also due to Paul Balmer for many helpful communications}. 
\section{Some Notations} 
First, we borrow some notations from (\cite{MSd}).

\begin{notations} \label{notaCate}{\rm
What follows would be our standard set up, throughout this
article. 

Throughout this paper,
$X:=(X,{\CO}_X)$ will denote a noetherian scheme, with $\dim X=d$.
We also assume $2$ is invertible in ${\CO}_X$.
Unless stated otherwise, we  assume that every coherent sheaf on    
$X$ is quotient of a locally free sheaf on $X$. 
This hypothesis is satisfied in the following two 
cases: (1) when $X$ is integral, locally factorial
and separated (see \cite[Ex. 6.8]{H}) and 
(2)  when $X$ has an ample invertible sheaf.  
To avoid  technicalities (see remark \ref{lowDepth}), we assume 
$$
d:=\dim X=\max\{depth({\CO}_{X,x}): x~{\rm is~a~closed~point} \}.
$$
This ensures that ${\CA}(X)$ has nonzero symmetric forms.  
We set up further notations:
\begin{enumerate}
\item $Coh(X)$ will denote the 
category of coherent ${\CO}_X-$modules and ${\SV}(X)$ will denote the 
full subcategory of locally free sheaves on $X$.
\item 
The full subcategory of objects ${\CF}$ in $Coh(X)$, so that ${\CF}$           
 has finite resolution by locally free sheaves, 
  will be denoted by ${\CB}:={\CB}(X)$.
  \item 
Also, let ${\CA}:={\CA}(X) \subseteq {\CB}(X)$
be   the full subcategory of objects  
${\CF}\in {\CB}(X)$ such that,${\CE}xt^i({\CF}, {\CO}_X)=0$
for all $i<d$. A coherent sheaf ${\CF}\in {\CA}(X)$ if and only if 
$co\dim(Supp({\CF}))=d$ and ${\CF}_x$ has finite
projective dimension for all closed points $x\in X$.  
%
   ({\it For our subsequent discussions, these two notations ${\CB}, {\CA}$
   will be of some importance.  In the absence of more standard
   notations or abbreviations of what they are, we use these two notations
   ${\CB}, ~{\CA}$. However, both ${\CB}$ and ${\CA}$   
   are  abelian subcategories 
   of $Coh(X)$,
   which justifies at least one of these two notations}.)

\item For any exact category ${\CC}$, $Ch^b({\CC}), D^b({\CC})$ will denote
the category of bounded chain complexes, and respectively, derived category.
If ${\CC}$ is a subcategory in an ambient abelian category ${\CC}'$ and ${\CH}$ is
a subcategory of ${\CC}'$, then $Ch^b_{{\CH}}({\CC})$ denotes the full subcategory of $Ch^b({\CC})$ consisting of complexes
with homologies in ${\CH}$. The derived category of 
$Ch^b_{{\CH}}({\CC})$ will be denoted by
$D^b_{{\CH}}({\CC})$, which is  
obtained by inverting quasi-isomorphisms in $Ch^b_{{\CH}}({\CC})$. 
However, $D^b_{{\CH}}({\CC})$  can also be viewed as the  
full subcategory of $D^b({\CC})$ consisting of objects from 
$Ch^b_{{\CH}}({\CC})$. Also, $K^b({\CC}), K^b_{{\CH}}({\CC})$ would 
denote the corresponding homotopy categories. 
For other similar notations, 
readers are referred to (\cite{W}). 
%
\item Denote objects ${\CE}_{\bullet}:=({\CE}_{\bullet}, \partial_{\bullet})
$ in $Ch^b(Coh(X))$  as:\\  
$\diagram
\cdots 0 \ar[r] & {\CE}_m \ar[r]^{\partial_m}
&{\CE}_{m-1} \ar[r] & \cdots  \ar[r]  &{\CE}_n \ar[r] & 0 \cdots
\enddiagram 
$ with $m > n$.
\item A  complex ${\CE}_{\bullet}$ is said to be supported on $[m,n]$ if
${\CE}_i=0$ unless $m\geq i\geq n$. 
\item For a complex ${\CE}_{\bullet}$ in $Ch^b({\SV}(X))$, 
${\CE}_{\bullet}^{\#}$ will denote the usual dual induced
by sheaf ${\CH}om(-,{\CO}_X)$ and $\varpi:{\CE}_{\bullet}\iso 
{\CE}_{\bullet}^{\#\#}$ will denote the
 evaluation map. 
\item Let $B_r=B_r({\CE}_{\bullet}):= \partial_{r+1}({\CE}_{r+1}) \subseteq {\CE}_r$
denote the module of  $r-$boundaries and ${\CZ}_r=
{\CZ}_r({\CE}_{\bullet}):= \ker(\partial_{r}) \subseteq {\CE}_r$
denote the module of $r-$cycles.
\item The $r^{th}-$homology of ${\CE}_{\bullet}$ will be 
denoted by ${\CH}_r:={\CH}_r({\CE}_{\bullet})$
and is defined by the exact sequence
$
\diagram
0\ar[r] & B_r({\CE}_{\bullet})\ar[r] & {\CZ}_r({\CE}_{\bullet})\ar[r] &
{\CH}_r\ar[r] & 0\\ 
\enddiagram
$.
\end{enumerate}
}
\end{notations}

\section{The functor by resolution}\label{zetaSection} 
Let $X$ be a noetherian scheme, as in (\ref{notaCate}). 
When $X=Spec(A)$ is affine, it is not difficult to
see that there is a natural functor
$\zeta: Coh(X) \lra K^{\geq 0}({\SV}(X))$, given by 
({\it choice of}) 
projective resolutions of finitely generated $A$-modules $M$ in $Coh(X)$.
If $M\in {\CB}$ has finite projective dimension then 
$\zeta(M)\in K^b({\SV}(X))$. 
When $X$
is not affine, morphisms $f:{\CF}\lra {\CG}$ in $Coh(X)$ 
do not naturally lift to  morphisms
of the corresponding ${\SV}(X)$-resolutions
of ${\CF}$ and ${\CG}$. 
 
However,
even when $X$ is not affine 
there are, in stead   functors  
$\zeta:{\CB}\lra D^b({\SV}(X))$, 
(resp. $\zeta:{\CB}(X)\lra D^{\geq 0}({\SV}(X))$)  
to the derived category, 
given  by ({\it the choices}) of  ${\SV}(X)$-resolutions
of objects in ${\CB}$ (resp. of $Coh(X)$). 
While this fact is well known among the experts,
there is no suitable reference,
to the best of our knowledge, that is accessible
to a wider community.
The purpose of this section is to give a proof of this fact 
for the benefit of such readers. 
We point out in \S ~\ref{resZetaSec},
that one can define and prove the existence of 
such a functor in the context of resolving subcategories of 
abelian categories. 

%
\bL \label{consTructRes}
Suppose $X$ is a  noetherian scheme as in (\ref{notaCate}).
Let $g:{\CF}\lra {\CG}$ be a morphism of coherent sheaves  
and 
$
\diagram 
\cdots \ar[r] & {\CG}_r \ar[r] & {\CG}_{r-1}\ar[r] &\cdots\ar[r] & {\CG}_0\ar[r] & {\CG} \ar[r] &0
\enddiagram 
$
be an exact sequence of 
coherent ${\CO}_X-$modules. 
Then, there is a resolution ${\CL}_{\bullet}$
of ${\CF}$ by objects in ${\SV}(X)$ 
and 
a morphism $g_{\bullet}:
{\CL}_{\bullet} \lra {\CG}_{\bullet}$ 
of complexes such that the diagram commutes:
$$
\diagram 
\cdots \ar[r] & {\CL}_r\ar[d]_{g_r} \ar[r] & {\CL}_{r-1}\ar[r] \ar[d]^{g_{r-1}}
&\cdots \ar[r]& {\CL}_0\ar[r]\ar[d]^{g_{0}} & {\CF}\ar[d]^g \ar[r] &0\\
\cdots\ar[r] & {\CG}_r \ar[r]_{\partial_r} & {\CG}_{r-1}\ar[r] &\cdots \ar[r]
& {\CG}_0\ar[r]_{\partial_0} & {\CG} \ar[r] &0\\ 
\enddiagram 
$$
Further, if ${\CF}$ is in ${\CB}(X)$ and ${\CG}_k=0~\forall k\gg 0$,
then ${\CL}_{\bullet}$ 
can be chosen to be in $Ch^b({\SV}(X))$. 
\eL 
\pf Consider the  pullback diagram
$$
\diagram 
\Gamma_0\ar[r]^{p_0}\ar[d]_{q_0} & {\CF}\ar[d]^g \\ 
{\CG}_0\ar@{->>}[r]_{\partial_0} & {\CG}.   
\enddiagram 
$$
Since $\partial_0$ is surjective, so is $p_0$. Since $\Gamma_0$ 
is coherent, by hypothesis, 
there is a 
surjective morphism $h_0:{\CL}_0\sur \Gamma_0$, where 
${\CL}_0 \in {\SV}(X)$. Let $d_0=p_0h_0$ and $g_0=q_0h_0$. 
So, we have a commutative diagram of exact sequences :
$$
\diagram 
0\ar[r] & {\CZ}_0\ar[r]\ar[d]_{G_0} &{\CL}_0\ar[r]^{d_0}\ar[d]^{g_0} 
& {\CF}\ar[d]^g \ar[r] & 0 \\ 
0\ar[r] & B_0\ar[r] &{\CG}_0\ar[r]_{\partial_0} & {\CG}\ar[r]
& 0 
\enddiagram 
$$
where $B_0=\ker(\partial_0)$, ${\CZ}_0=\ker(d_0)$, and $G_0$ 
is the restriction of $g_0$.  
We extend the diagram as follows:
$$
\diagram 
\Gamma_1\ar[r]^{p_1}\ar[d]_{q_1} 
& {\CZ}_0\ar[r]\ar[d]_{G_0} &{\CL}_0\ar[r]^{d_0}\ar[d]^{g_0} 
& {\CF}\ar[d]^g \ar[r] & 0 \\ 
{\CG}_1\ar@{->>}[r]_{\partial_1}
& B_0\ar[r] &{\CG}_0\ar@{->>}[r]_{\partial_0} & {\CG}\ar[r] & 0 
\enddiagram 
$$
Here $\Gamma_1$ is the pullback.  
Note $p_1$ is surjective.
Now, the proof of the first part is complete by iteration of this process. 
The latter statement follows because the process terminates, in that case.  
$\eop$ 
The following is another version of (\ref{consTructRes}) that will be
of some interest subsequently.
\bL\label{quasiIsoInPX}
Suppose $X$ is a  noetherian scheme, as in (\ref{notaCate}).
Let ${\CG}_{\bullet}$ be a complex in $Ch^+(Coh(X))$.
Then there is a quasi-isomorphism $g_{\bullet} :{\CL}_{\bullet}\lra 
{\CG}_{\bullet}$,
with 
 ${\CL}_{\bullet}$ in $Ch^{+}({\SV}(X))$.  
Further, if ${\CG}_{\bullet}$ is in $Ch^b(Coh(X))$,
such that all the homologies ${\CH}_i({\CG}_{\bullet}) \in {\CB}(X)$.
Then 
${\CL}_{\bullet}$ can be chosen in
$Ch^{b}({\SV}(X))$.
\eL 
\pf ({\it With}
\S ~\ref{resZetaSec} {\it in mind,
we avoid local argument
and stick to formal arguments}). We write ${\CG}_{\bullet}$ as 
$\diagram \cdots \ar[r] & 
{\CG}_{n+1} \ar[r] & {\CG}_{n} \ar[r]^{\partial_n} & {\CG}_{n-1} 
\ar[r] & \cdots \enddiagram$. We assume that ${\CG}_i=0~\forall i<0$.
Denote $B_i:=image(\partial_{i+1})$ and
$Z_i:=\ker(\partial_i)$.
%
First, by hypothesis, there is a surjective map $g_0: {\CL}_0\sur 
{\CG}_0$, with ${\CL}_0\in {\SV}(X)$. 
 With $d_0=\partial_0 g_0$ and ${\CZ}_0=\ker(d_0)$, we have 
 a commutative diagram of exact sequences:
$$
\diagram
0\ar[r]&{\CZ}_0\ar@{->>}[d]^{f_0}\ar[r]& 
{\CL}_0\ar[r]^{d_0}\ar@{->>}[d]^{g_0}& {\CH}_0({\CG}_{\bullet})\ar@{=}[d]\ar[r] & 0\\
0\ar[r]&Z_0\ar[r]& {\CG}_0\ar[r]_{\partial_0}& {\CH}_0({\CG}_{\bullet})\ar[r] &0\\ 
\enddiagram
\quad f_0~{\rm where~is~the~restriction~of}~g_0.  
$$
It follows $f_0$ is surjective. Inductively,
we choose surjective maps $g_n: {\CL}_n\sur {\CG}_n$ and differentials
$d_n:{\CL}_n\lra {\CL}_{n-1}$,
such that (1) ${\CL}_n\in {\SV}(X)$ and (2) with ${\CZ}_n:=\ker(d_n)$,
the restriction map 
$f_n:{\CZ}_n \sur Z_n$  of $g_n$ to ${\CZ}_n:=\ker(d_n)$  
is surjective. So, we have the following
commutative diagram of exact sequences: 
$$
\diagram
0\ar[r]&{\CZ}_n\ar@{->>}[d]^{f_n}\ar[r]& 
{\CL}_n\ar[r]^{d_n}\ar[d]^{g_n}& {\CL}_{n-1}\ar[d]^{g_{n-1}}\\
0\ar[r]&Z_n\ar[r]& {\CG}_n\ar[r]_{\partial_n}& {\CG}_n.\\ 
\enddiagram
$$
Consider the pullback diagram: 
$$
\diagram
0\ar[r]&{\CZ}_{n+1}\ar@{->>}[d]^{f_{n+1}}\ar[r]
&{\CL}_{n+1}\ar@{-->>}[d]^{\varphi}\ar[rrrd]^{d_{n+1}}&&&\\
0\ar[r] & Z_{n+1}\ar[r]\ar@{=}[d]& 
\Gamma_{n+1}\ar@{->>}[d]_{q}\ar@{->>}[r]^{p} 
& f_n^{-1}(B_n)\ar@{->>}[d]^{f_n'}\ar@{^(->}[r] 
& {\CZ}_n\ar@{^(->}[r]\ar@{->>}[d]^{f_n} &{\CL}_n
\ar[d]^{g_n}\\
0\ar[r] & Z_{n+1}\ar[r]&
{\CG}_{n+1}\ar@{->>}[r]_{\partial_{n+1}} & B_n\ar@{^(->}[r]&Z_n\ar@{^(->}[r] &{\CG}_n\\  
\enddiagram
$$
In this diagram (1) 
$\Gamma_{n+1}$ is the pullback of $(\partial_{n+1}, f_n')$
and (2) we chose a surjective map 
$\varphi: {\CL}_{n+1}\sur \Gamma_{n+1}$ with ${\CL}_{n+1}\in {\SV}(X)$.
We have the following:
(1)  By the properties of pull back 
$p, q$ are surjective and $\ker(p)\cong Z_{n+1}$.  
(2) Let $g_n=q \varphi$ and the differential $d_{n+1}$ is defined as shown. 
Denote 
${\CZ}_{n+1}:=\ker(d_{n+1})$.
(3)  Since $\varphi, p$ are surjective, $Image(d_{n+1})= f_n^{-1}(B_n)$.
(4) Since $f_n$ is  surjective,
the homology map ${\CH}(g_n): {\CH}_n({\CL}_{\bullet})
\sur {\CH}_n({\CG}_{\bullet})$ is surjective.
It also follows that ${\CH}(g_n)$ is injective.  
Therefore, ${\CH}(g_n): 
{\CH}_n({\CL}_{\bullet})
\iso  {\CH}_n({\CG}_{\bullet})$ is an isomorphism.
(5) By Snake lemma, it also follows
that $f_{n+1}$ is surjective. 
This completes the proof of the first part of the lemma.  

Finally, if all the homologies ${\CH}_i({\CG}_{\bullet})\in {\CB}(X)$,
since ${\CB}(X)$ is an exact
catgory and every epimorphism in ${\CB}(X)$ is admissible, 
 ${\CZ}_i({\CL}_{\bullet}), B_i({\CL}_{\bullet})\in {\CB}(X)$. 
The left tail of ${\CL}_{\bullet}$ becomes a resolution of
$B_n({\CL}_{\bullet})$ for some $n\gg 0$. Since  ${\CG}_i=0~\forall ~i\gg 0$, 
we can truncate ${\CL}_{\bullet}$.
\pic $\eop$


Now, we proceed to establish that
given a morphism $g:{\CF}\lra {\CG}$
of coherent sheaves, and resolutions ${\CL}_{\bullet}, {\CE}_{\bullet}$,
respectively, of ${\CF},{\CG}$, $g$  
lifts to a morphism ${\CL}_{\bullet}\lra {\CE}_{\bullet}$ in the derived 
category.
\bL \label{pullBres}
Suppose $X$ is a noetherian scheme, as in (\ref{notaCate}). 
Let $g:{\CF}\lra {\CG}$ be a 
morphism of coherent sheaves and ${\CF}_{\bullet},~{\CG}_{\bullet}$
be complexes  in $Ch^{\geq 0}({Coh}(X))$ such that 
$$ 
{\CH}_0({\CF}_{\bullet})
= {\CF},\quad {\CH}_0({\CG}_{\bullet})
= {\CG},~~~ {\rm and}~~~ 
{\CH}_i({\CF}_{\bullet})={\CH}_i({\CG}_{\bullet})=0~\forall~i\neq 0,
$$
Consider ${\CG}$ as a complex in $Ch^{\geq 0}({Coh}(X))$ concentrated at degree
zero and let $\Gamma_{\bullet}$ be the pullback:
$$
\diagram
\Gamma_{\bullet}\ar[r]^t\ar[d]_G & {\CF}_{\bullet} \ar[d]\\
{\CG}_{\bullet} \ar[r]& {\CG.}\\ 
\enddiagram \quad {\rm Then,~\exists ~a~commutative~diagram}
\quad 
\diagram
{\CF}_{\bullet}\ar@{->>}[d]
& \Gamma_{\bullet} \ar[l]_t\ar[r]^G\ar@{->>}[d] 
& {\CG}_{\bullet}\ar@{->>}[d] \\
{\CF}\ar@{=}[r] &{\CF}\ar[r]_g & {\CG}\\ 
\enddiagram
$$
such that $G, t$ are morphisms in $Ch^{\geq 0}({Coh}(X))$
and $t:\Gamma_{\bullet}\lra {\CF}_{\bullet}$ is a quasi-isomorphism.
\eL
\pf We establish the notations as in the diagram:
$$
\diagram 
\cdots \ar[r] & {\CF}_r \ar[r]^{d_r} & {\CF}_{r-1}\ar[r]^{d_{r-1}} 
&\cdots \ar[r]& {\CF}_0\ar[r]^{d_0} & {\CF}\ar[d]^g \ar[r] &0\\
\cdots\ar[r] & {\CG}_r \ar[r]_{\partial_r} & {\CG}_{r-1}\ar[r] &\cdots \ar[r]
& {\CG}_0\ar[r]_{\partial_0} & {\CG} \ar[r] &0\\ 
\enddiagram 
$$
where the horizontal lines are exact. 
Recall that the pullback of complexes is constructed by taking 
degree wise pullback. 
So, $\Gamma_i=0~\forall~i<0$. At degree $i=0$, of  the pullback diagram
is given by:
$$
\diagram  
\Gamma_0\ar[r]^{t_0}\ar[d]_{G_0} & {\CF}_0\ar[d]^{gd_0}\\
{\CG}_0\ar@{->>}[r]_{\partial_0} & {\CG}\\
\enddiagram 
\qquad \qquad {\rm Since}~\partial_0~{\rm is~surjective,~so~is}~t_0. 
$$
Therefore,
$\partial'_0:=gt_0: \Gamma_0\sur {\CF}$ is surjective. 
We have $\ker(\partial'_0)=\ker(d_0)\oplus \ker(\partial_0)$.
However, $\forall~i>0~~\Gamma_i={\CF}_i\oplus {\CG}_i$
and the differentials $\partial'_i=d_i\oplus 
\partial_i:\Gamma_i\lra \Gamma_{i-1}$.
\pic $\eop$

\bC \label{2ResRiso}
Suppose $X$ is a noetherian scheme, as in (\ref{notaCate}).
Let ${\CF}$ be  an object
in ${Coh}(X)$  
and ${\CF}_{\bullet}, {\CF}_{\bullet}'$
be two resolutions of 
${\CF}$ in the chain complex category 
$Ch^{\geq 0}(X)$.
%
Then ${\CF}_{\bullet}\cong {\CF}_{\bullet}'$ in 
the derived category $D^{+}(Coh(X))$. 
\eC
\pf It follows immediately, by an application of 
corollary \ref{pullBres}, with $g=1_{\CF}$, since in this case
$G$ is also a quasi-isomorphism.
$\eop$ 


\bC \label{functMap}
Suppose $X$ is a noetherian scheme,
as in (\ref{notaCate}). 
Under the hypotheses and notations in lemma \ref{pullBres}, there is 
a complex ${\CL}_{\bullet}\in Ch^{\geq 0}({\SV}(X))$,
a quasi-isomorphism $\tau:{\CL}_{\bullet}\lra {\CF}_{\bullet}$,
and a
morphism $\gamma:{\CL}_{\bullet}\lra {\SV}_{\bullet}$ 
such that the diagram 
$$
\diagram
{\CF}_{\bullet}\ar@{->>}[d]& {\CL}_{\bullet} \ar[l]_{\tau}\ar[r]^{\gamma}
\ar@{->>}[d] 
& {\CG}_{\bullet}\ar@{->>}[d] \\
{\CF}\ar@{=}[r] &{\CF}\ar[r]^g & {\CG}\\ 
\enddiagram
\qquad \quad {\rm commutes}. 
$$
Further, if ${\CF}\in {\CB}(X)$ and ${\CG}_{\bullet} \in Ch^b({Coh}(X))$, 
then ${\CL}_{\bullet}$ can be chosen to be in $Ch^b({\SV}(X))$.
\eC 
\pf 
We use all the notations as in lemma \ref{pullBres}. By lemma 
\ref{consTructRes}, applied to the identity morphism $Id:{\CF}\lra {\CF}$
and $\Gamma_{\bullet}$, there
a complex ${\CL}_{\bullet}\in Ch^+({\SV}(X))$ and a
quasi-isomorphism $\epsilon:{\CL}\lra \Gamma_{\bullet}$.
Now, the corollary is established,
with the choices $\tau=t\epsilon$ and $\gamma=G\epsilon$. 
The last statement also follows from the corresponding statement
in lemma \ref{consTructRes}.
\pic $\eop$
For  morphisms $g:{\CF}\lra {\CG}$ in $Coh(X)$, 
a "lift"  in the derived category is defined as follows. 
\bD \label{defFunt(g)}
Suppose $X$ is a noetherian scheme, in (\ref{notaCate}).
Let $g:{\CF}\lra {\CG}$ be a
morphism of coherent sheaves on $X$. Suppose  ${\CE}_{\bullet}\sur {\CF},
{\CQ}_{\bullet}\sur {\CG}$
are two resolutions,  
where ${\CE}_{\bullet}, {\CQ}_{\bullet}$ are in $Ch^{\geq 0}({\SV}(X))$.
By corollary \ref{functMap}, 
there is
a complex ${\CL}_{\bullet}\in Ch^{\geq 0}({\SV}(X))$,
a quasi-isomorphism $\tau:{\CL}_{\bullet}\lra {\CE}_{\bullet}$,
and a
morphism $\gamma:{\CL}_{\bullet}\lra {\CQ}_{\bullet}$
such that the diagram
$$
\diagram
{\CE}_{\bullet}\ar@{->>}[d]& {\CL}_{\bullet} \ar[l]_{\tau}\ar[r]^{\gamma}
\ar@{->>}[d] 
& {\CQ}_{\bullet}\ar@{->>}[d] \\
{\CF}\ar@{=}[r] &{\CF}\ar[r]^g & {\CG}\\ 
\enddiagram
\quad {\rm commutes. ~Define}~~\zeta(g):{\CE}_{\bullet} \lra {\CQ}_{\bullet}
\quad {\rm by}
\quad \zeta(g)=\gamma \tau^{-1}
$$
as a morphism $ Mor_{D^+({\SV}(X))}\left({\CE}_{\bullet},
{\CQ}_{\bullet}\right)$ in the derived category.
%
Further, if ${\CE}_{\bullet}, {\CQ}_{\bullet}\in Ch^b({\SV}(X))$ 
then $\zeta(g)$ is defined in $D^b({\SV}(X))$.
\eD 
We now establish that $\zeta(g)$ is well defined.
\bL \label{wdnessOfMap}  
Under the set up of definition \ref{defFunt(g)}, 
$\zeta(g)$ is well defined in $D^+({\SV}(X))$. 
Further, $\zeta(g)$ is in $D^b({\SV}(X))$,
if  ${\CE}_{\bullet},{\CQ}_{\bullet} \in Ch^b({\SV}(X))$. 
\eL 
\pf 
As in (\ref{pullBres}), consider ${\CG}$
as a complex in $Ch^+({Coh}(X))$ and let $\Gamma_{\bullet}$ be the pullback:
$$
\diagram
\Gamma_{\bullet}\ar[r]^t\ar[d]_G & {\CE}_{\bullet} \ar[d]\\
{\CQ}_{\bullet} \ar[r]& {\CG}.\\ 
\enddiagram
\quad\qquad {\rm Then,~the~diagram}\qquad
\diagram
{\CE}_{\bullet}\ar@{->>}[d]
& \Gamma_{\bullet} \ar[l]_t\ar[r]^G\ar@{->>}[d] 
& {\CQ}_{\bullet}\ar@{->>}[d] \\
{\CF}\ar@{=}[r] &{\CF}\ar[r]_g & {\CG}\\ 
\enddiagram \quad {\rm commutes}
$$
and $t$ is a quasi-isomorphism. 
Given two pairs $({\CL}_{\bullet}, \tau, \gamma)$ and
$({\CL}'_{\bullet}, \tau', \gamma')$ as in definition \ref{defFunt(g)},
from the properties of the pullback it follows, that there  are  
maps $\epsilon:{\CL}_{\bullet}\lra \Gamma_{\bullet}$, 
$\epsilon':{\CL}_{\bullet}'\lra \Gamma_{\bullet}$ such that the diagrams 
$$
\diagram
&&{\CE}_{\bullet}\ar[rd] & \\
{\CL}_{\bullet} \ar@{-->}[r]^{\epsilon}\ar[rru]^{\tau}\ar[rrd]_{\gamma} 
& \Gamma_{\bullet}\ar[ru]_t\ar[rd]^G  && {\CG} \\ 
&&{\CQ}_{\bullet}  \ar[ru]&\\
\enddiagram 
\quad {\rm and}\quad 
\diagram
&&{\CE}_{\bullet}\ar[rd] & \\
{\CL}_{\bullet}' \ar@{-->}[r]^{\epsilon'}\ar[rru]^{\tau'}\ar[rrd]_{\gamma'} 
& \Gamma_{\bullet}\ar[ru]_t\ar[rd]^G  && {\CG} \\ 
&&{\CQ}_{\bullet}  \ar[ru]&\\
\enddiagram 
\quad {\rm commute.} 
$$
Note that $\epsilon, \epsilon'$ 
are quasi-isomorphisms. 
By Ore condition (see proof of \cite[10.4.1]{W}) in $Coh^+(X)$ there is 
$\Delta_{\bullet}\in Ch^+(Coh(X))$ and quasi-isomorphisms   
$\mu: \Delta_{\bullet}\lra {\CL}_{\bullet}$ and 
$\mu': \Delta_{\bullet}\lra {\CL}_{\bullet}'$ such that 
$\epsilon \mu=\epsilon'\mu'$.   
%
By 
(\ref{quasiIsoInPX}), there is quasi-isomorphism $\nu:{\CU}_{\bullet} \lra
\Delta_{\bullet}$, with ${\CU}_{\bullet}\in Ch^+({\SV}(X))$.
$$
{\rm With}~~\eta=\mu\nu,~\eta'=\mu'\nu,
\quad {\rm ~ ~ the~diagram~}\quad 
\diagram 
&&{\CL}_{\bullet} \ar[rd]^{\epsilon} \\ 
{\CU}_{\bullet} \ar[r]^{\nu}
\ar[rru]^{\eta} \ar[rrd]_{\eta'}
&\Delta_{\bullet}\ar[ru]_{\mu}\ar[rd]^{\mu'}&& \Gamma_{\bullet}\\
&&{\CL}_{\bullet}' \ar[ru]_{\epsilon'} \\ 
\enddiagram
\quad {\rm commutes}. 
$$
All the morphisms in this latter diagram
are quasi-isomorphisms.
It follows, 
$$
\tau\eta 
=\tau'\eta',
~~ 
\gamma \eta 
= \gamma'\eta'. 
\quad {\rm Therefore, ~the~diagram}~~
\diagram 
&{\CL}_{\bullet}\ar[dr]^{\gamma}\ar[ld]_{\tau}&\\
{\CE}_{\bullet}&{\CU}_{\bullet}\ar[u]^{\eta}\ar[d]_{\eta'}\ar[r]^{\gamma \eta}
\ar[l]_{\tau\eta}
&{\CQ}_{\bullet}\\ 
&{\CL}'_{\bullet}\ar[ru]_{\gamma'}\ar[lu]^{\tau'}&\\
\enddiagram 
$$
commutes. 
So, $\gamma \tau^{-1} =\gamma'(\tau')^{-1} $. This establishes that  $\zeta(g)$ is well defined in $D^+({\SV}(X))$.
Similarly, it follows that  $\zeta(g)$ is well defined in $D^b({\SV}(X))$,
if ${\CE}_{\bullet}, {\CQ}_{\bullet}$ are in $Ch^b({\SV}(X))$. 
\pic $\eop$

In deed, as in the affine case, $\zeta(g)$ is the unique lift of
$g$, in the following sense. 
\bT\label{anytaugamma}
Consider the set up as in definition \ref{defFunt(g)}.
Then, $\zeta(g): {\CE}_{\bullet}\lra {\CQ}_{\bullet}$ is the unique
morphism $\eta:{\CE}_{\bullet}\lra {\CQ}_{\bullet}$ such that ${\CH}_0(\eta)=g$.
\eT 
\pf Suppose $\eta$ is as in the statement of the corollary. 
Then, $\eta=\gamma'(\tau')^{-1}$, where $\tau', \gamma'$ is given by 
the commutative diagram: 
$$
\diagram
{\CE}_{\bullet}\ar@{->>}[d]& {\CL}_{\bullet}' \ar[l]_{\tau'}\ar[r]^{\gamma'}
& {\CQ}_{\bullet}\ar@{->>}[d] \\
{\CF}\ar[rr]_g & & {\CG}\\ 
\enddiagram
\quad {\rm where}~~{\CL}_{\bullet}' \in Ch^+({\SV}(X)),~ {\rm and}~~\tau'~~
{\rm is ~a~ quasi-isomorphism}.
$$
 By replacing sheaves at the negative degrees of ${\CL}_{\bullet}'$ by zero
and 
${\CL}'_0$ by $\ker({\CL}'_0\lra {\CL}'_{-1})$, we obtain a   
complex ${\CL}_{\bullet}$ and a quasi-isomorphim $\epsilon: {\CL}_{\bullet}
\lra {\CL}_{\bullet}'$,
such that (1) ${\CL}_{\bullet}$ fits in the diagram of 
the definition \ref{defFunt(g)}, with $\tau=\tau'\epsilon$ and 
$\gamma=\gamma'\epsilon$. 
It follows,  $\zeta(g)=\gamma\tau^{-1}= \gamma'(\tau')^{-1}$.  \pic $\eop$

Now, $\zeta$ behaves in a functorial manner,
for the morphisms, as follows. 
\bP \label{funtOfZeta}
Let $X$ be a noetherian scheme as in (\ref{notaCate}). 
Let $g_0:{\CF}_0\lra {\CF}_1$ 
and $g_1:{\CF}_1\lra {\CF}_2$ be 
morphisms of coherent sheaves on $X$. For $i=0, 1,2$
let ${\CE}_{\bullet}^i\sur {\CF}_i$ be resolutions of 
${\CF}_i$, with ${\CE}_{\bullet}^i$ in $Ch^{\geq 0}({\SV}(X))$
{\rm (resp. in $Ch^{b}({\SV}(X))$)}.
%
Then $\zeta(g_1g_0)=\zeta(g_1)\zeta(g_0)$, where $\zeta(g_0),
\zeta(g_1), \zeta(g_1g_0)$ 
are as defined in (\ref{defFunt(g)}).
 \eP
\pf Consider the following commutative diagram  
$$
\diagram
&&{\CE}_{\bullet}\ar[dr]^{\gamma}\ar[dl]_{\tau}&&\\ 
{\CE}^0_{\bullet}\ar@{->>}[d]&{\CL}^0_{\bullet}\ar[l]_{\tau_0}
\ar[r]^{\gamma_0}\ar@{->>}[d]
&{\CE}^1_{\bullet}\ar@{->>}[d]&{\CL}^1_{\bullet}\ar@{->>}[d]\ar[l]_{\tau_1}
\ar[r]^{\gamma_1}
&{\CE}^2_{\bullet}\ar@{->>}[d]\\ 
{\CF}_0\ar@{=}[r] & {\CF}_0\ar[r]_{g_0} & {\CF}_1\ar@{=}[r] 
& {\CF}_1 \ar[r]_{g_1} & {\CF}_2\\  
\enddiagram ~~~~~~~~~~~~~~~.
$$
In this diagram, $({\CL}^i_{\bullet}, \tau_i,\gamma_i)$ are obtained by
corollary \ref{functMap}, for $i=0,1$, corresponding to $g_i$,
where ${\CL}^i_{\bullet}$ are in $Ch^{\geq 0}({\SV}(X))$ and $\tau_i$
are quasi-isomorphism. Further, $({\CE}, \tau,\gamma)$ 
is obtained by application of the Ore condition in $Ch^+({\SV}(X))$ 
(resp. in $Ch^b({\SV}(X))$, where $\tau$ is a quasi-isomorphism.
By (\ref{anytaugamma}), it follows,
$\zeta(g_1g_0)=$
$$
(\gamma_1\gamma)(\tau_0\tau)^{-1}=
(\gamma_1\tau_1^{-1})(\tau_1\gamma)(\tau_0\tau)^{-1}
=\zeta(g_1)(\gamma_0\tau)(\tau_0\tau)^{-1}=\zeta(g_1)\gamma_0\tau_0^{-1}
=\zeta(g_1)\zeta(g_0). 
$$
\pic $\eop$


We are finally ready to give a
complete definition of  the functors $\zeta:{Coh}(X)\lra D^+({\SV}(X))$,
and $\zeta:{\CB}(X)\lra D^b({\SV}(X))$, 
as follows. 

\bD \label{finalDefZeta} 
Let $X$ be a noetherian scheme, as in (\ref{notaCate}).
Define 
$\zeta:{\CB}(X)\lra D^b({\SV}(X))$ as follows 
\bE
\item For each object ${\CF}\in {\CB}(X)$, by {\rm (\ref{consTructRes})}, 
choose a resolution ${\CE}_{\bullet}\sur {\CF}$
where ${\CE}_{\bullet}$
is in $Ch^{b}({\SV}((X)))$.   
Define
$$
\zeta({\CF}):={\CE}_{\bullet} \quad in
\quad D^b({\SV}(X)).
$$
{\rm For simplicity, we make a convention that if the length of a minimal 
${\SV}(X)$-resolution of 
${\CF}$ is $r$, then by choice $\zeta({\CF})_i=0 ~\forall ~i>r$}.
\item For  morphisms $g:{\CF}\lra {\CG}$ in ${\CB}(X)$,
define $\zeta(f):\zeta({\CF})\lra \zeta({\CG})$ 
as in  \ref{defFunt(g)}.
%
\eE
The functor, $\zeta:Coh(X)\lra D^+({\SV}(X))$ is defined similarly. 
\eD 
It is immediate from proposition \ref{funtOfZeta} that the
functors $\zeta$ defined above are well defined.  


\subsection{The Resolving category and fomalism} \label{resZetaSec} 
Much of the arguments above can be formulated in the realm  
of resolving subcategories of
abelian categories as follows. Before we proceed,
we define resolving subcategories of abelian categories.
\bD\label{resCatDef}{\rm 
Suppose ${\CC}$ is an abelian category. An exact subcategory
${\SV}$ of ${\CC}$ is called a {\bf resolving subcategory}  
if,
(1) ${\SV}$ is closed under
direct summand and direct sum, 
(2) every epimorphism ${\SV}$ is admissible, and
(3) given any object $M$ in ${\CC}$ there is an epimorphism 
${\CE}\sur {\CG}$, for some ${\CE}\in {\SV}$. 
}\eD 

\begin{remark}\label{resremark}{\rm
Suppose ${\SV}$ is a resolving subcategory of an abelian category ${\CC}$.
Following are some immediate and obvious comments:

\noindent (1) Any object $M$ is ${\CC}$ has a resolution by objects in ${\SV}$.
(2)  The category of chain complexes $Ch^b({\SV})$, $Ch^+({\SV})$, 
$Ch({\SV})$
of objects in ${\SV}$
and the   homotopy categories $K^b({\SV})$, $K^+({\SV})$, and $K({\SV})$ 
are defined, as usual. 
(3)  The derived categories 
$D^b({\SV})$, $D^+({\SV})$, and $D({\SV})$
are defined
by inverting the quasi-isomorphisms in the corresponding
homotopy category ({\it without any regard to any other structure}). 
Consult the proof of \cite[10.4.1]{W} 
on localization, which makes sense in our context.  
}
\end{remark} 
\bP \label{formalZeta}
{\rm
Suppose $\CC$ is an abelian category and ${\SV}$ is a 
resolving  subcategory.
$$
{\rm Let}\quad {\CB}:={\CB}({\SV}):=
\{{\CF}\in {\SV}: {\CF} {\rm ~has~a~finite~{\SV}-resolution}\}
\quad {\rm be~the~full~subcategory.} 
$$
Assume, given an object ${\CF}\in {\CB}$ and a
${\SV}$-resolution ${\CE}_{\bullet}\sur {\CF}$, 
the cycle objects  
$Z_n:=\ker({\CE}_{n}\lra {\CE}_{n-1})\in {\SV}$~~$\forall~n\gg 0$.
\bE
\item\label{1OfForZeta}
${\CB}$ is an exact subcategory and every epimorphism in 
${\CB}$ is admissible. 
\item \label{2OfForZeta}
There is a natural  functor
$$
\zeta: {\CB}({\SV}) \lra D^b(\SV) \qquad 
{\rm 
as~defined~ in~(\ref{defFunt(g)})}. 
$$
Likewise, there is a functor $\zeta: {\CC} \lra D^+(\CP)$ (without any condition on $Z_n$).
\item \label{3OfForZeta} All of \S ~\ref{zetaSection}
remain valid, if we replace $Coh(X)$ by
${\CC}$, ${\SV}(X)$ by ${\SV}$ and ${\CB}(X)$ by ${\CB}(\SV)$.  
\eE 
}\eP
\pf The statements  (\ref{2OfForZeta}, \ref{3OfForZeta}) follows
as in the cases of noetherian schemes. 
We need to give a proof of (\ref{1OfForZeta}).
Let
$
\diagram 0\ar[r]& K \ar[r]_f & M \ar[r]_g & C\ar[r] & 0\\
\enddiagram
$ be an exact sequence in ${\CC}$. Inductively, a resolution (possibly infinite)
 of this sequence in $Ch^{\ge 0}({\SV})$, as follows.  
Consider the diagram
$$
\diagram
0\ar[r]& L_0\ar@{->>}[d]_{\partial_0} \ar[r] & P_0\ar@{->>}[d]_{\varphi_0}
\ar[r] & Q_0\ar@{=}[d]\ar[r] & 0\\ 
0\ar[r]& K\ar@{=}[d] \ar[r] & \Gamma_0\ar@{->>}[d]_{d_0'}
\ar[r]_{\gamma_0}
& Q_0\ar@{->>}[d]^{d_0}\ar[r] & 0\\ 
 0\ar[r]& K \ar[r]_f & M \ar[r]_g & C\ar[r] & 0\\
\enddiagram
$$
In this diagram, $d_0:Q_0\sur C$ is a surjective morphism,
 is the pullback of $(d_0,g)$
and $\varphi: P_0\sur \Gamma_0$ is a surjective morphism,
with $Q_0, P_0\in {\SV}$. Also, $L_0:=\ker(\gamma_0\varphi_0) 
\in {\SV}$. By Snake lemma $\partial_0$ is 
surjective. Write $\delta_0=d'_0\varphi_0$. Again,
by Snake lemma\\  
$\diagram 
0\ar[r]& \ker(\partial_0)  \ar[r] & \ker(\delta_0) \ar[r] 
& \ker(d_0)\ar[r] & 0 \\
\enddiagram $
is exact and the process continues. 

Now, if $K, C\in {\CB}({\SV})$ then, by hypothesis, the resolutions 
$L_{\bullet}\sur K, Q_{\bullet}\sur C$ terminate, hence
process stops and  
 $P_{\bullet}\in Ch^b({\SV})$. So, ${\CB}({\SV})$ is an exact subcategory.
Similarly, every epimorphims in ${\CB}({\SV})$ is admissible.
\pic $\eop$  
\section{Duality Properties}\label{secWittA(X)} 
In this section $X$ will denote a noetherian scheme,
as in (\ref{notaCate}), with $\dim X=d$.  
Also, ${\CA}(X)$ will denote the abelian subcategory of $Coh(X)$,
as define in (\ref{notaCate}). 
With a view on section \ref{FormalDiv} on formalism, we take the 
formal approach to the proofs, as opposed to local. 
First, we define a duality on ${\CA}(X)$.
\bD \label{extChoice}
Suppose $X$ is a noetherian scheme with $\dim X=d$,  as in 
(\ref{notaCate}). 
We fix a choice of an injective resolution ${\SI}_{\bullet}$ of ${\CO}_X$.
For a coherent sheaf ${\CF}$, and integers $n\geq 0$,
we define ${\CE}xt^n({\CF}, {\CO}_X):={\CH}^n({\CF}, {\SI}_{\bullet})$. 
\eD

\bL\label{dualCheck} 
Suppose $X$ is a noetherian scheme with $\dim X=d$
and ${\CA}(X)$ is as in (\ref{notaCate}). 
Define a functor 
$^{\vee}:{\CA} \lra {\CA}$
 by ${\CF}^{\vee}:= {\CE}xt^d\left({\CF}, {\CO}_X\right)$.
Then, $^{\vee}$ defines a duality on ${\CA}$. 
\eL
\pf We need to establish that, for objects ${\CF}\in {\CA}(X)$,  there is 
a natural isomorphism $\tilde{\varpi}: {\CF}\iso {\CF}^{\vee\vee}$. 
Let ${\CE}_{\bullet}:=\zeta({\CF})$
 as defined in 
(\ref{finalDefZeta}).
We denote this complex as 
$$
\diagram 
0\ar[r] & {\CE}_d \ar[r]^{\partial_d} & \cdots \ar[r]& {\CE}_1\ar[r]^{\partial_1}  & {\CE}_0\ar[r]^{\partial_0} &  {\CF} \ar[r] & 0\\
\enddiagram 
$$
So, there is a natural isomorphism $co\ker(\partial_d^*)\iso
{\CF}^{\vee}$ (see \cite{H}). 
Since ${\CF}\in {\CA}(X)$ we have ${\CE}xt^i({\CF}, {\CO}_X)=0$ for 
all $i\neq d$. 
Therefore, the dual ${\CE}^{\#}$ yields a resolution:  
$$
\diagram 
0\ar[r] & {\CE}_0^* \ar[r] & \cdots \ar[r]& {\CE}_{d-1}^*\ar[r]  & {\CE}_d^*\ar[r] &  
{\CF}^{\vee}\ar[r] & 0
\enddiagram
$$
Dualizing again, we have the following diagram of exact sequences
$$
\diagram
0\ar[r] & {\CE}_d \ar[r]\ar[d]_{ev} & \cdots \ar[r]
& {\CE}_{1}\ar[r]\ar[d]_{ev}  & {\CE}_0\ar[r]\ar[d]_{ev} &  
{\CF}\ar[r]\ar@{-->}[d]^{\varpi_0}_{\wr} & 0\\ 
0\ar[r] & {\CE}_d^{**} \ar[r] & \cdots \ar[r]& {\CE}_{1}^{**}\ar[r]  & 
{\CE}_0^{**}\ar[r] &  
co\ker(\partial_0^{**})\ar[r]
& 0\\ 
\enddiagram 
$$ 
The isomorphism $\varpi_0$ is induced by the evaluation maps. 
There is also a natural isomorphism $\omega_1:co\ker(\partial_0^{**})\iso 
{\CF}^{\vee\vee}$. Hence, $\tilde{\varpi}:= \varpi_1\varpi_0:
{\CF}\iso {\CF}^{\vee\vee}$ is a natural isomorphism.  
\pic $\eop$


The derived category $D^b({\SV}(X)))$ has a triangulated structure,
with the duality induced by the duality  $-^*:={\CH}om(-,{\CO}_X)$,
which we will denote by $\#$.
We have particular interest in the subcategory $D^b_{{\CA}}({\SV}(X)))$. 
It was pointed out in (\cite{MSd}) that $D^b_{{\CA}}({\SV}(X)))$ may not
be closed under cone construction, but it was stable under duality. This
works for any noetherian scheme, as follows. 
\bL \label{Ext12Hmt}
Suppose $X$ is a noetherian scheme, with $\dim X=d$ {\rm (the condition in (\ref{notaCate}) that every coherent
sheaf on $X$ is quotient of a sheaf in ${\SV}(X)$, is not needed.)}
Also, ${\CA}:={\CA}(X)$ be as in (\ref{notaCate}). 
Let ${\CE}_{\bullet}$ be a complex in $Ch^b_{{\CA}}({\SV}(X))$. 
Then, for all integers $r\in{\BZ}$, 
there are natural isomorphism
\begin{equation} \label{Ebat2Hx}
{\CE}xt^i\left(\frac{{\CE}_r}{B_r}, {\CO}_X \right)\iso
\left\{\begin{array}{ll}
{\CE}xt^d({\CH}_{r+i-d}({\CE}_{\bullet}), {\CO}_X) &1\leq i\leq d\\0&i\geq d+1\\\end{array}
\right. 
\end{equation}
where $B_r\subseteq {\CE}_r$ denote the boundary, as in (\ref{notaCate}).
\eL
\pf We assume ${\CE}_r=0~\forall r<0$. 
If ${\CH}_0({\CE}_{\bullet})=0$ then there is nothing to prove. 
So, assume the ${\CH}_0({\CE}_{\bullet})\neq 0$. Then, 
$\frac{{\CE}_0}{B_0}={\CH}_0({\CE}_{\bullet})$ and equation 
(\ref{Ebat2Hx}) holds.
Now assume that  equation 
(\ref{Ebat2Hx}) holds for degree $r$.
We  prove it for degree 
$r+1$.
We have two exact sequences
$$
\diagram
0\ar[r] & B_r \ar[r] & {\CE}_r \ar[r] & \frac{{\CE}_r}{B_r}\ar[r] & 0
&{\rm and}& 
0\ar[r] & {\CH}_{r+1}
\ar[r] &  \frac{{\CE}_{r+1}}{B_{r+1}} \ar[r] & B_r\ar[r] & 0
\enddiagram
$$
The long exact sequence of the first exact sequence yields the 
following isomorphisms: 
$$
\forall~i\geq 1\quad {\CE}xt^{i}(B_r, {\CO}_X) \iso 
{\CE}xt^{i+1}\left( \frac{{\CE}_r}{B_r}, {\CO}_X \right) 
\iso
\left\{\begin{array}{ll}
{\CE}xt^d({\CH}_{r+(i+1)-d}, {\CO}_X) &1\leq i\leq d-1\\0&i\geq d\\\end{array}
\right. 
$$
Whether ${\CH}_{r+1}=0$ or ${\CH}_{r+1}\neq d$, 
${\CE}xt^{i}({\CH}_{r+1}, {\CO}_X)=0~\forall~i\neq d$. 
It follows from the second exact sequence 
$$
{\CE}xt^{i}\left( \frac{{\CE}_{r+1}}{B_{r+1}}, {\CO}_X \right)\iso
\left\{
\begin{array}{ll}{\CE}xt^i(B_r,{\CO}_X)&if~i\leq d-1\\
{\CE}xt^d({\CH}_{r+1},{\CO}_X)&if~i=d\\0&if~i\geq d+1 \end{array}
\right.
$$
This establishes the lemma. $\eop$ 
 
\bT \label{dual4ChbPX}
Let $X$ be a noetherian scheme as in (\ref{Ext12Hmt}) and  
${\CE}_{\bullet}$ be a complex in $Ch^b_{{\CA}}({\SV}(X))$. Then,
for $r\in{\mathbb Z}$, there is a canonical isomorphism 
$
\eta_{{\CE}_{\bullet}}: 
{\CH}_{-r}({\CE}_{\bullet}^{\#}) 
\iso {\CH}_{r-d}({\CE}_{\bullet})^{\vee}. 
$
In particular, 
$Ch^b_{{\CA}}({{\SV}}(X))$ is stable under duality.
Further, $\eta_{{\CE}_{\bullet}}$ is natural with
respect to morphisms  $f:{\CE}_{\bullet}\lra {\CE}_{\bullet}'$  in 
$Ch^b_{{\CA}}({\SV}(X))$.   
\eT 
\pf 
First, we have the following commutative diagram of exact sequences:
{\scalefont{0.7}
$$
\diagram
0\ar[r] & \left(\frac{{\CE}_{r-1}}{B_{r-1}}\right)^* \ar[r] \ar@{=}[d]
& \left({\CE}_{r-1}\right)^* \ar[r]\ar@{=}[d]
& (B_{r-1})^* \ar[d]^{\wr}\ar[r]
& {\CE}xt^1\left(\frac{{\CE}_{r-1}}{B_{r-1}} ,{\CO}_X\right)\ar@{-->}[d]_{\wr}^{\psi}\ar[r] & 0\\
0\ar[r] & \left(\frac{{\CE}_{r-1}}{B_{r-1}}\right)^* \ar[r] 
& \left({\CE}_{r-1}\right)^* \ar[r]
& \left(\frac{{\CE}_{r}}{B_{r}}\right)^* \ar[r]
& H_{-r}({\CE}_{\bullet}^*) \ar[r] & 0
\enddiagram
$$
}
The isomorphism $\psi$ is induced. Now, the first part of the theorem 
follows by composing $\psi$ with the isomorphism, as follows: 
${\CE}xt^1\left(\frac{{\CE}_{r-1}}{B_{r-1}} ,{\CO}_X\right)\iso 
{\CE}xt^d({\CH}_{r-d},{\CO}_X)
$
$$
\diagram
{\CE}xt^1\left(\frac{{\CE}_{r-1}}{B_{r-1}} ,{\CO}_X\right)\ar[d]^{\psi}
\ar[r]^{\sim} &  
{\CE}xt^d({\CH}_{r-d},{\CO}_X)\ar[dl]^{\eta_{{\CE}_{\bullet}}} \\
 H_{-r}({\CE}_{\bullet}^{\#})&\\ 
\enddiagram
$$
The latter part follows because all the isomorphisms are natural.
\pic $\eop$  
\section{The Main Isomorphism Theorem} \label{mainIsoSEC} 
Having established the that $D^b_{\CA}({\SV}(X))$ is stable under duality,
we discuss Witt group of $D^b_{\CA}({\SV}(X))$. 
We set up some basic framework, analogous to (\cite{MSd}). 
\bE
\item For a complex ${\CE}_{\bullet}$ in $Ch^b({\SV}(X))$, 
$T_u({\CE}_{\bullet})$ will denote the unsigned translation,
and $T_s({\CE}_{\bullet})$ the standard translation which changes the sign of
the differential. 
\item  Denote 
the shifted Derived categories:
$
T^nD^b_{{\CA}}({\SV}(X))^{\pm}_u:= \left(D^b_{{\CA}}({\SV}(X)),T^n_uo\#, 1,
\pm \varpi \right) 
$
and 
$
T^nD^b_{{\CA}}({\SV}(X))^{\pm}_s:= \left(D^b_{{\CA}}({\SV}(X)),T^n_so\#, 1,
\pm \varpi \right)
$
where $\varpi$ is the evaluation map. 
\item
As was pointed out, $D^b_{\CA}({\SV}(X))$ may not have
a triangulated structure. 
Readers are also 
referred to (\cite{MSd}), for a definition of Witt groups of 
subcategories of triangulated categories with duality.
Note that this definition is very similar to  
that of Witt groups of exact catgories in \cite{TWGII}, requiring 
the lagrangians to be admissible monomorphisms. 
As usual, we denote $W^n(D^b_{{\CA}}({\SV}(X))^{\pm}_u):=
W(T^nD^b_{{\CA}}({\SV}(X))^{\pm}_u)$ and 
$W^n(D^b_{{\CA}}({\SV}(X))^{\pm}_s):=
W(T^nD^b_{{\CA}}({\SV}(X))^{\pm}_s)$. 
\eE
Now, we prove that the functor $\zeta:{\CA}(X)\lra D^b_{{\CA}}({\SV}(X))$
induces isomorphisms of 
Witt groups. 
First, we prove $\zeta$ induces a homomorphism
of Witt groups.

\bT \label{zWDWittmap}
Suppose $X$ is a noetherian scheme, as in (\ref{notaCate}), with 
$d=dim X$. Let ${\CA}:={\CA}(X)$ be as in (\ref{notaCate}).
Then, the functor $\zeta:{\CA}\lra D^b_{{\CA}}({\SV}(X))$ induces 
a well defined homomorphism
$$
W(\zeta): W({\CA}, ^{\vee}, \pm\tilde{\varpi}) \lra 
W^d\left(D^b_{{\CA}}({\SV}(X))^{\pm}_u\right)
$$
\eT
\pf We will only give the proof for the $+$-duality.
We first define the homomorphism $W(\zeta):W({\CA}, ^{\vee}, \tilde{\varpi}) \lra 
W^d\left(D^b_{{\CA}}({\SV}(X))^+_u)\right)$ and then
prove that it is well defined.
Suppose $({\CF},\varphi_0)$ is a symmetric form in 
$({\CA}, ^{\vee}, \tilde{\varpi})$. 
Write ${\CX}_{\bullet}:=
\zeta({\CF})$. By choice, ${\CX}_i=0~\forall~i>d$.  
 Then, ${\CX}_{\bullet}^{\#}$ gives a resolution of ${\CF}^{\vee}$.
It is also easy to see that $\zeta$ preserves the duality on ${\CA}$.
By the uniqueness (\ref{anytaugamma})  of the lifts of morphisms in 
${\CA}(X)$, the identity $\varphi_0=\varphi_0^{\vee}\tilde{\varpi}$ produces 
a symmetric from $\zeta(\varphi_0):{\CX}_{\bullet}\iso {\CX}_{\bullet}^{\#}$.
Therefore, the association $({\CF}, \varphi_0)\mapsto
({\CX}_{\bullet}, \zeta(\varphi_0))$  
induces a homomorphism
 $MW(\zeta): MW\left({\CA}, ^{\vee}, \tilde{\varpi}\right)\lra 
MW(T^d\left(D^b_{{\CA}}({\SV}(X))^+_u)\right))$ of monoids of symmetric spaces. 

Now we want to prove that $MW(\zeta)$ maps neutral spaces to neutral 
spaces. 
Suppose
$({\CF},\varphi_0)$  is a neutral space in $\left({\CA}, ^{\vee}, \tilde{\varpi}\right)$.
So, there is an exact sequence (i.e. a lagrangian)
$\diagram
0\ar[r] & {\CG} \ar[r]^{\alpha_0} & {\CF} \ar[r]^{\alpha_0^{\vee}\varphi_0}
& {\CG}^{\vee} \ar[r] &0. 
\enddiagram $
As above, write ${\CX}_{\bullet}:=\zeta({\CF})$. We would like to 
show $({\CX}_{\bullet}, \zeta(\varphi_0))$ is neutral. 
It will be convenient to work with forms
without  denominators. 
 By (\ref{functMap} or \ref{defFunt(g)}), $\zeta(\varphi_0)
 =\gamma \tau^{-1}: {\CX}_{\bullet}\iso {\CX}_{\bullet}^{\#}$, where $\tau, \gamma$ are quasi-isomorphisms in 
 $Ch^b({\SV}(X)$, ${\CE}_{\bullet}\in Ch^+({\SV}(X))$, 
 such that
the diagram 
\begin{equation} \label{zWDEqn} 
\diagram
{\CX}_{\bullet} \ar@{->>}[d] 
& {\CE}_{\bullet}\ar@{->>}[d]\ar[l]_{\tau}\ar[r]^{\gamma}
& {\CX}_{\bullet}^{\#} \ar@{->>}[d] \\
{\CF}\ar@{=}[r] & {\CF} \ar[r]_{\varphi_0} & {\CF}^{\vee}\\ 
\enddiagram
\quad{\rm commutes.~~With}~~\varphi=\tau^{\#}\gamma~~{\rm the~diagram}\quad 
\diagram
{\CE}_{\bullet}\ar@{->>}[r]\ar[d]_{\varphi}
& {\CF} \ar[d]^{\varphi_0}\\
{\CE}_{\bullet}^{\#}\ar@{->>}[r]& {\CF}^{\vee} \\
\enddiagram 
\end{equation}
commutes and 
 $\varphi$ is a quasi-isomorphism in $Ch^b({\SV}(X))$,
without denominator. 
Clearly, $({\CE}_{\bullet}, \varphi)$ and $({\CX}_{\bullet}, \zeta(\varphi_0))$
are isometric. 

Now that $MW(\zeta)({\CF}, \varphi_0)$ is given by 
a denominator free quasi-isomorphism $\varphi:{\CE}_{\bullet}\lra
{\CE}_{\bullet}^{\#}$, the rest of the argument is borrowed from 
(\cite{DWG}, {\it for extra details see} \cite{MSd}). 
We outline the proof to point out the subtleties. 
Let ${\CL}_{\bullet}:=\zeta({\CG})$. By choice,
${\CL}_0=0$ unless $d\geq i\geq 0$. 
Then, ${\CL}_{\bullet}^{\#}$ yields a resolution of 
${\CG}^{\vee}$. 
Since $\varphi=\zeta(\varphi_0)$, by (\ref{funtOfZeta}), 
the above exact sequence lifts 
to an exact sequence 
$
\diagram
0\ar[r] & {\CL}_{\bullet} \ar[r]^{\alpha} & {\CE}_{\bullet}
\ar[r]^{\alpha^{\#}\varphi}
& {\CL}_{\bullet}^{\#} \ar[r] &0  
\enddiagram 
$. 
Since the homotopy category $K^b({\SV}(X))$ has a
triangulated structure, 
we  can embed $\alpha$ in an exact triangle in $T^dK^b({\SV}(X))$
and 
 obtain the following 
diagram and a morphism $s$, as follows:
$$
\diagram
{\CL}_{\bullet} \ar[r]^{\alpha}\ar[d] 
& {\CE}_{\bullet} \ar[r]^j\ar[d]^{\alpha^{\#}\varphi}  
& {\CV}_{\bullet}\ar@{-->}[d]^s\ar[r]^k & T({\CL}_{\bullet})\ar[d]^0\\
0\ar[r] & {\CL}_{\bullet}^{\#}\ar@{=}[r]
& {\CL}_{\bullet}^{\#}\ar[r] & 0.\\ 
\enddiagram
\qquad \qquad {\rm in}\quad T^dK^b({\SV}(X)). 
$$
From the long exact sequence of homologies, it follows 
${\CH}_0({\CV}_{\bullet})\cong {\CG}^{\vee}$ and 
${\CH}_i({\CV}_{\bullet})=0~\forall~i\neq 0$. So, 
$s$ is an isomorphism in $D^b_{{\CA}}({\SV}(X))$.
Therefore, we get the exact triangle  
$$
\diagram
{\CL}_{\bullet} \ar[r]^{\alpha} & {\CE}_{\bullet} 
\ar[r]^{\alpha^{\#}\varphi}
& {\CL}_{\bullet}^{\#}\ar[r]^{ks^{-1}} & T({\CL}_{\bullet}) \\ 
\enddiagram 
\quad {\rm in} \quad D^b_{{\CA}}({\SV}(X)). 
$$
Translating this triangle, and with $w=-T^{-1}(ks^{-1})$,
we get 
the exact triangle
$$
\diagram
T^{-1}{\CL}_{\bullet}^{\#}\ar[r]^w
& {\CL}_{\bullet} \ar[r]^{\alpha} & {\CE}_{\bullet} 
\ar[r]^{\alpha^{\#}\varphi}
& {\CL}_{\bullet}^{\#}\\
\enddiagram 
\quad {\rm in} \quad T^dD^b_{{\CA}}({\SV}(X)).  
$$
It remains to show that $T^{-1}w^{\#}=w$. 
Routine analysis
shows that this is equivalent to showing
$T(k^{\#})s=T(s^{\#})k$. Indeed, $s:{\SV}_n={\CE}_n\oplus {\CL}_{n-1}
\lra {\CL}^*_{d-n}$ is a morphism of complexes. Therefore all the maps in 
this equation $T(k^{\#})s=T(s^{\#})k$ are morphisms of complexes. 
It would suffice to show  that this identity holds in $T^dK^b({\SV}(X))$,
which can be done exactly as in (\cite{DWG, MSd}).
This establishes that $\zeta$ induces a homomorphism of the Witt
groups. $\eop$ 

Following lemma will be helpful for subsequent discussion.
\bL \label{sdtSupp} Suppose $X$ is a noetherian scheme 
with $\dim X=d$.   
Suppose $({\CE}_{\bullet}, \varphi)$
is a symmetric from 
in $T^dD^b_{{\CA}}({\SV}(X))$. Assume, for some $n\geq 0$,  
${\CH}_{-n}({\CE}_{\bullet})\neq 0$ and
${\CH}_{i}({\CE}_{\bullet})=0~\forall~i<-n$.
Then, $({\CE}_{\bullet}, \varphi)$ is 
isometric,  to a symmetric from 
$({\CE}_{\bullet}', \varphi')$ such that ${\CE}'_i=0$, unless
$n+d\geq i\geq -n$. 
\eL
\pf ({\it We give a proof that applies to} \S  ~\ref{FormalDiv}.)
We use the notations from (\ref{notaCate}).
Note the cycle sheaf
${\CZ}_{-n}({\CE}_{\bullet})\in Ch^b({\SV}(X))$. Replacing,
${\CE}_{i}$ by zero, when $i<-n$ and 
${\CE}_{-n}$ by ${\CZ}_{-n}({\CE}_{\bullet})$,  
we obtain a complex ${\CE}'_{\bullet}$, with ${\CE}'_i=0
~ \forall ~i<-n$ and a 
quasi-isomorphism  to $\eta: {\CE}_{\bullet}'\lra 
{\CE}_{\bullet}$.
Therefore, replacing $({\CE}_{\bullet}, \varphi)$
by $({\CE}_{\bullet}',  \eta^{\#}\varphi \eta)$, we assume that 
${\CE}_i=0~\forall ~i<-n$. 
  
By the duality theorem \ref{dual4ChbPX}, $\forall~k>n$,  we have
${\CH}_k({\CE}_{\bullet}) \cong {\CH}_k({\CE}_{\bullet}^{\#})\cong
{\CH}_{-k}({\CE}_{\bullet})^{\vee}=0$. 
Inductively, we show that we can assume ${\CE}_i=0$ for 
all $i\geq n+d+1$. We will assume ${\CE}_i=0$ for 
all $i\geq n+d+2$. Write $K:=\ker(\partial_{n+d+1}^*)$.
Then, $K\in {\SV}(X)$. Dualizing the  exact sequences
$
\diagram
0\ar[r] & K \ar[r] & {\CE}_{n+d}^*
\ar[r]^{\partial_{n+d+1}^{*}} & {\CE}_{n+d+1}^*\ar[r] &0 
\enddiagram
$,
we get the exact sequence 
$
\diagram
0\ar[r]  &  {\CE}_{n+d+1}\ar[r]^{\partial_{n+d+1}} & {\CE}_{n+d}\ar[r]
& K^*\ar[r] &0.  
\enddiagram $ 
Since $K^*\in {\CV}(X)$, by replacing ${\CE}_i$ by zero for 
$i\geq n+d+1$ and ${\CE}_{n+d}$ by $K^*$, we get a new complex
${\CE}_{\bullet}'$ and a quasi-isomorphism $\eta:{\CE}_{\bullet}\lra 
{\CE}_{\bullet}'$. 
It follows $({\CE}_{\bullet}, \varphi)$ is isometric to 
$({\CE}_{\bullet}',  (\eta^{-1})^{\#}\varphi \eta^{-1})$.
\pic $\eop$

%


As in (\cite{DWG, TWGII}),
the proof of surjectivity of $W(\zeta)$ is done by
reduction of support of the symmetric form by application 
of the sublagrangian theorem (\cite[4.20]{TWGI}).
In the non-affine case, both the construction 
and the proof that it is a sublagrangian
 require further finesse. The following lemma will be useful.  
%
%
\bL \label{itIsSubL}
Suppose $X$ is as in (\ref{notaCate}). Let ${\CL}_{\bullet},
{\CG}_{\bullet}\in Ch^b({\SV}(X))$ be complexes
and \\ $\eta_{\bullet}:{\CL}_{\bullet}
\lra {\CG}_{\bullet}$ be a morphism, as in the  diagram
$$
\diagram 
\cdots \ar[r] &
0 \ar[r] & {\CL}_{n}\ar[d]_{\eta} \ar[r] & {\CL}_{n-1}\ar[r] \ar[d]^{\eta}
&\cdots \ar[r]& {\CL}_0\ar[r]\ar[d]^{\eta_{0}} 
& 0 \ar[r] &\cdots\\
\cdots \ar[r] &
0\ar[r] & {\CG}_{n} \ar[r] & {\CG}_{n-1}\ar[r] &\cdots \ar[r]
& {\CG}_0\ar[r] & {\CG}_{-1} \ar[r] &\cdots\\ 
\enddiagram 
$$
such that
\bE
\item ${\CH}_r({\CG}_{\bullet})=0~\forall r\geq 0$.
\item ${\CH}_r({\CL}_{\bullet})=0~\forall r\neq 0$ and ${\CL}_r=0~\forall~r<0$.
\eE
Then, $\eta=0$ in $D^b({\SV}(X))$. 
\eL 
\pf We use the notations in (\ref{notaCate}). 
Imitating the arguments in (\ref{wdnessOfMap}), let $\Gamma_{\bullet}$
be the complex:
$
\diagram
\cdots \ar[r] & {\CL}_2\oplus {\CG}_2 \ar[r] & {\CL}_1\oplus {\CG}_1
 \ar[r] & {\CL}_0\oplus {\CZ}_0({\CG}_{\bullet})\ar[r] 
 & 0\ar[r] & 0\ar[r] &\cdots 
\enddiagram
$
Let $t:\Gamma_{\bullet}\lra {\CL}_{\bullet}$ and 
$g:\Gamma_{\bullet}\lra {\CG}_{\bullet}$  
be the projection maps. Then, $t$ is a quasi-isomorphism.
We have the following commutative diagram of morphisms:
$$
\diagram
& {\CL}_{\bullet}\ar[rd]^{\eta}\ar[ld]_{1}\ar[d]^{(1,\eta)}& \\
{\CL}_{\bullet} & \Gamma_{\bullet}\ar[l]_t\ar[r]^g & {\CG}_{\bullet}\\ 
& {\CL}_{\bullet}\ar[ru]_{0}\ar[lu]^{1}\ar[u]_{(1,0)}& \\ 
\enddiagram
$$
Since, $t$ is quasi-isomorphism, so is $(1,\eta)$. Rest of the 
arguments in (\ref{wdnessOfMap}) works, which uses Ore condition,
together with (\ref{consTructRes}). Note $\Gamma_{\bullet},
{\CL}_{\bullet}$ are resolutions, so good enough to apply 
arguments in (\ref{wdnessOfMap}). 
\pic $\eop$


Now we prove that  $W(\zeta)$ is surjective, as follows.
\bP \label{surjWzeta}
Let $X$ be  a noetherian scheme, with $\dim X=d$,
as in (\ref{notaCate}). 
Then, the homomorphism
$
W(\zeta): W({\CA}, ^{\vee}, \pm\tilde{\varpi}) \lra 
W^d\left(D^b_{{\CA}}({\SV}(X))^{\pm}_u)\right)
$
is surjective. 
\eP 
\pf We point out the subtleties involved here in the
non-affine case, beyond the arguments in (\cite{MSd}).
We will only consider the case of  $+$duality 
and the case of skew duality follows similarly. 
As usual, the proof is done by 
reducing the length of the width of the forms.
Suppose $x=[({\CE}_{\bullet}, \varphi)]\in 
W^d\left(D^b_{{\CA}}({\SV}(X))_u)\right)$, represented by the 
a symmetric form $({\CE}_{\bullet}, \varphi)$. First,
$\varphi=fs^{-1}$ for some quasi-isomorphism $s:{\CE}_{\bullet},\lra 
{\CE}_{\bullet}$ for some ${\CE}_{\bullet}' \in Ch^b({\SV}(X))$. 
Replacing  $({\CE}_{\bullet}, \varphi)$ by  
$({\CE}_{\bullet}', s^{\#}\varphi s)$,
we assume that $\varphi: {\CE}_{\bullet}\lra {\CE}_{\bullet}^{\#}$
is a quasi-isomorphism (without denominator) in $Ch^b({\SV}(X))$.

Assume, for some $n>0$, suppose ${\CH}_{-n}({\CE}_{\bullet})\neq 0$ and 
${\CH}_{i}({\CE}_{\bullet})=0~\forall~i<-n$. 
We will prove that there is a symmetric form $({\CR}_{\bullet},\psi)$
such that $x=[({\CE}_{\bullet}, \varphi)]= [({\CR}_{\bullet},\psi)]$ 
and ${\CH}_i({\CR}_{\bullet})=0$ unless $n-1\geq i\geq -(n-1)$.
By lemma \ref{sdtSupp}, we assume that ${\CE}_i=0$ unless $n+d\geq i
\geq -n$. Further, by duality theorem \ref{dual4ChbPX},
${\CH}_i({\CE}_{\bullet})\cong 
{\CH}_i({\CE}_{\bullet}^{\#})
\cong {\CH}_{-i}({\CE}_{\bullet})^{\vee}=0$ for all $i>d$.
So, the left tail ${\CE}_{\bullet}$
is a resolution of $\frac{{\CE}_n}{B_n({\CE}_{\bullet})}$.
By lemma \ref{consTructRes}, the morphism ${\CH}_n({\CE}_{\bullet}) 
\hra \frac{{\CE}_n}{B_n({\CE}_{\bullet})}$
induces morphism of complexes, as follows:
$$
\diagram 
0 \ar[r] & {\CL}_{n+d}\ar[d]_{\nu} \ar[r] & {\CL}_{n+d-1}\ar[r] \ar[d]^{\nu}
&\cdots \ar[r]& {\CL}_n\ar[r]\ar[d]^{\nu_{n}} 
& {\CH}_n({\CE}_{\bullet})\ar@{^(->}[d] \ar[r] &0\\
0\ar[r] & {\CE}_{n+d} \ar[r] & {\CE}_{n+d-1}\ar[r] &\cdots \ar[r]
& {\CE}_n\ar[r] & \frac{{\CE}_n}{B_n({\CE}_{\bullet})} \ar[r] &0\\ 
\enddiagram 
$$
where ${\CL}_i\in {\SV}(X)$ and both the lines are exact. 
Since $\nu_n$ maps to ${\CZ}_{n}({\CE}_{\bullet}))$,
with ${\CL}_i=0~\forall~i\leq n-1$, it extends to 
morphism $\nu:{\CL}_{\bullet}\lra {\CE}_{\bullet}$ 
in $Ch^b({\SV}(X))$. We claim $\nu$ is a sublagrangian
of $({\CE}_{\bullet},\varphi)$. To see this, write $\eta=\nu^{\#}\varphi
\nu:{\CL}_{\bullet}\lra {\CL}_{\bullet}^{\#}$. We can write $\eta$ as
follows
$$
\diagram 
0 \ar[r] & {\CL}_{n+d}\ar[d]_{\eta} \ar[r] & {\CL}_{n+d-1}\ar[r] \ar[d]^{\eta}
&\cdots \ar[r]& {\CL}_n\ar[r]\ar[d]^{\eta_{n}} 
& 0 \ar[d]\ar[r] &\cdots\\
0\ar[r] & {\CL}_{-n}^{\#} \ar[r] & {\CL}_{-(n-1)}^{\#}\ar[r] &\cdots \ar[r]
& {\CL}_{d-n}^{\#}\ar[r] & {\CL}_{d-n+1}^{\#} \ar[r] &\cdots\\ 
\enddiagram 
$$
Since, only at degree $-n$, $L_{\bullet}^{\#}$ has a non zero homology,
the second line is exact at degrees $i\geq n$. 
By lemma \ref{itIsSubL}, $\eta=0$ in $D^b({\SV}(X))$.
Rest of the arguments in \cite{MSd} works, which we outline
briefly for completeness. 
%
As in \cite[4.20]{TWGI}, embed $\nu$ in an exact triangle  
as in the first line in the diagram (\ref{usualDiagam}) and 
complete the commutative diagram
\begin{equation} \label{usualDiagam}
\xymatrix{
T^{-1} {\CN}_{\bullet} \ar[r]^{\nu_0}\ar[d]_{T^{-1}\mu_0^{\#}}
& {\CL}_{\bullet} \ar[r]^{\nu}\ar[d]^{\mu_0}
& {\CE}_{\bullet}\ar[d]^{\varphi} \ar[r]^{\nu_1} 
& {\CN}_{\bullet}\ar[d]^{\mu_0^{\#}}\\
T^{-1} {\CL}_{\bullet}^{\#} \ar[r] & {\CN}_{\bullet}^{\#}\ar[d]^{\mu_1} \ar[r]
& {\CE}_{\bullet}^{\#} \ar[r] & {\CL}_{\bullet}^{\#}\\
& {\CR}_{\bullet}\ar[d]^{\mu_2}  &  & \\
& T{\CL}_{\bullet}  &  & \\
} 
\end{equation}
where the second line the dual of the first line, by choice $\mu_0$ is
"very good" (see \cite{TWGI} for definition and existence) and 
${\CR}_{\bullet}$ is the cone of $\mu_0$. We consider this as 
a diagram in $D^b({\SV}(X))$, which is a triangulated category. 
By \cite[4.20]{TWGI} there is a symmetric form $\psi:{\CR}_{\bullet}
\iso {\CR}_{\bullet}^{\#}$ such that $({\CR}_{\bullet},-\psi)$ is 
Witt equivalent to $({\CE}_{\bullet},\varphi)$ in $D^b({\SV}(X))$.
This is shown by exhibiting a lagrangian 
$\lambda: 
{\CN}_{\bullet}^{\#} \lra ({\CE}_{\bullet}, -\varphi)\perp 
({\CR}_{\bullet}, \psi)$. We would have to show that ${\CR}_{\bullet}$ 
is in $D^b_{{\CA}}({\SV}(X))$ and $\lambda$ is a lagrangian in 
$D^b_{{\CA}}({\SV}(X))$. It suffices to show, 
\bE
\item ${\CN}_{\bullet}$,
${\CN}_{\bullet}^{\#}$ and ${\CR}_{\bullet}$ are in $D^b_{{\CA}}({\SV}(X))$,
\item and 
%
${\CH}_i({\CR}_{\bullet})=0$ unless $n-1\geq i\geq -(n-1)$.
\eE
These are established by  writing down the long exact sequences
of homologies 
of the three exact triangles in 
the diagram (\ref{usualDiagam})  
and  using the fact that 
${\CL}_{\bullet}, {\CL}_{\bullet}^{\#}$ have only one nonzero homology.
({\it With a view to section \ref{FormalDiv} on formalism, we 
avoid local argument.})
By lemma \ref{sdtSupp}, we can further assume that ${\CR}_i=0$
unless $(n-1)+d\geq i\geq -(n-1)$. 
Using induction there is a symmetric from $({\CQ}_{\bullet}, \omega)$
such that $x=[({\CE}_{\bullet}, \varphi)]=[({\CQ}_{\bullet}, \omega)]$,
with ${\CQ}_i=0$
unless $d\geq i\geq 0$. By (\ref{dual4ChbPX}),
${\CQ}_{\bullet}$ is a resolution 
of ${\CH}_0({\CQ}_{\bullet})$. Further, $\omega$ induces 
a symmetric from $\omega_0:{\CH}_0({\CQ}_{\bullet})\iso {\CH}_0({\CQ}_{\bullet})^{\vee}$. By definition,
$W(\zeta)([({\CH}_0({\CQ}_{\bullet},\omega_0)]=[({\CQ}_{\bullet}, \omega)]
=x$. So, $W(\zeta)$ is surjective. 
\pic
$\eop$

\vspace{2mm}
The following is the main d\'{e}vissage theorem for noetherian 
schemes.
\bT \label{divScheme}
Let $X$ be  a noetherian scheme, with $\dim X=d$,
as in (\ref{notaCate}). 
Then, the homomorphism 
$
W(\zeta): W({\CA}(X), ^{\vee}, \pm\tilde{\varpi}) \lra 
W^d\left(D^b_{{\CA}}({\SV}(X))^{\pm}_u\right)
$
defined in (\ref{zWDWittmap}), 
is an isomorphism. 
\eT
\pf It follows from (\ref{surjWzeta}) that 
$W(\zeta)$ is surjective. 
The proof of injectivity is similar to that in the affine case (\cite{MSd}),
by an application of (\cite[4.1]{TWGII}). We outline the
proof briefly, in the case of $+$duality. 
Suppose $({\CF},\varphi_0)$ is a symmetric form in
$\left({\mathcal A}, ^{\vee}, \tilde{\varpi} \right)$ and
$W(\zeta)([({\CF},\varphi_0)])=0$. 
By definition \ref{zWDWittmap} (see diagram \ref{zWDEqn})
there is a finite resolution 
${\CE}_{\bullet}\sur {\CF}$, with ${\CE}_{\bullet}\in Ch^b({\SV}(X))$
and a quasi-isomorphism (without denominator) 
$\varphi:{\CE}_{\bullet}\lra {\CE}_{\bullet}^{\#}$,
which lifts $\varphi_0$ and $W(\zeta)([({\CF},\varphi_0)])
=[({\CE}_{\bullet}, \varphi)]=0$.
So, $({\CE}_{\bullet}, \varphi)$ is neutral in 
$D^b_{{\CA}}\left({\SV}(X)\right)$. Going through the same 
argument in (\cite[5.12]{MSd}) 
there is a hyperbolic from $({\CQ},\varphi_1)\perp ({\CQ},-\varphi_1)$
such that $({\CE},\varphi)\perp ({\CQ},\varphi_1)\perp ({\CQ},-\varphi_1)$ 
is neutral in $D^b_{{\CA}}\left({\SV}(X)\right)$. 
Therefore, it follows 
$$
({\CU}_{\bullet}, \beta):= 
\left({\CE}_{\bullet} \oplus {\CQ}_{\bullet}\oplus {\CQ}_{\bullet}^{\#}, 
\left(\begin{array}{ccc}\varphi_0 & 0 & 0\\ 0& 0 & 1 \\
0& 1 &0 \end{array} \right)\right) 
\quad \text{is neutral in} \quad \left(D^b_{{\mathcal A}}\left({\SV}(X) \right)^{\pm}_u\right).
$$
Therefore, 
$(U_{\bullet}, \varphi)$ has a lagrangian $({\CL}_{\bullet}, \alpha)$
given by the following exact triangle
$$
\diagram
T^{-1}{\CL}_{\bullet}^{\#} \ar[r]^w & {\CL}_{\bullet}  \ar[r]^{\alpha}  
& {\CU}_{\bullet}  \ar[r]^{\alpha^{\#}\varphi} 
& {\CL}_{\bullet}^{\#}  
\enddiagram 
\qquad {\rm with} \quad T^{-1}w^{\#}=w.  
$$
By the duality theorem \ref{dual4ChbPX}, ${\CH}_{-r}({\CU}_{\bullet})\iso 
{\CH}_{d-r}({\CU}_{\bullet})^{\vee}$  
and ${\CH}_{-r}({\CL}_{\bullet})\iso 
{\CH}_{d-r}({\CL}_{\bullet})^{\vee}$, for all $r\in {\BZ}$.   
With these identification, the exact sequence of the homologies
of 
the exact triangle reduces to 
{\scalefont{0.7}
$$\xymatrixcolsep{5pc}\xymatrix{
\ar[r] & {\CH}_{-2}({\CL}_{\bullet})^{\vee} \ar[r]^{{\CH}_1(w)} 
& {\CH}_{1}({\CL}_{\bullet}) \ar[r]^{{\CH}_1(\alpha)} 
& {\CH}_1({\CU}_{\bullet}) \ar[r]^{{\CH}_{-1}(\alpha)^{\vee} \smallcirc{0.5} {\CH}_1(\beta)} 
& {\CH}_{-1}({\CL}_{\bullet})^{\vee} \\
}$$
$$\xymatrixcolsep{5pc}\xymatrix{
\ar[r]^{{\CH}_0(w)} &  {\CH}_0({\CL}_{\bullet}) \ar[r]^{{\CH}_0(\alpha)} &
{\CH}_0({\CU}_{\bullet}) \ar[r]_{{\CH}_{0}(\alpha)^{\vee} \smallcirc{0.5} {\CH}_0(\beta)} &
{\CH}_0({\CL}_{\bullet})^{\vee} \ar[r]_{{\CH}_0(w)^{\vee}} & {\CH}_{-1}({\CL}_{\bullet})^{\vee \vee} \\
}$$
$$\xymatrixcolsep{5pc}\xymatrix{
\ar[r]_{{\CH}_1(\beta)^{\vee} \smallcirc{0.5} {\CH}_{-1}(\alpha)^{\vee \vee}} &
{\CH}_1({\CU}_{\bullet})^{\vee} \ar[r]_{{\CH}_1(\alpha)^{\vee}} &
{\CH}_1({\CL}_{\bullet})^{\vee} \ar[r]_{{\CH}_1(w)^{\vee}} &
{\CH}_{-2}({\CL}_{\bullet})^{\vee \vee} \ar[r]_{{\CH}_2(\beta)^{\vee} \smallcirc{0.5} 
{\CH}_{-2}(\alpha)^{\vee \vee}} & \\
}
$$
}
This exact sequence is "symmetric" and  hence \cite[4.1]{TWGII} applies.
Since the sequence is exact, 
it follows 
$
[({\CH}_0({\CU}_{\bullet}), {\CH}_0(\varphi))]=[(0,0)]= 0$
in $W({\CA}(X), ^{\vee}, \tilde{\varpi})$. 
However,
$$
({\CH}_0({\CU}_{\bullet}), {\CH}_0(\beta)) =
\left({\CH}_0({\CE}_{\bullet}) \oplus {\CH}_0({\CQ}_{\bullet})\oplus {\CH}_0({\CQ}_{\bullet}^{\#}), 
\left(\begin{array}{ccc}{\CH}_0(\varphi_0) & 0 & 0\\ 0& 0 & 1 \\
0& 1 &0 \end{array} \right)\right) 
$$
$$
= ({\CF}, \varphi_0)
\perp \left({\CH}_0({\CQ}_{\bullet})\oplus {\CH}_0({\CQ}_{\bullet})^{\vee}, 
\left(\begin{array}{cc} 0 & 1 \\ 1 &0 \end{array} \right)\right)  
$$
So, we have
$$
[({\CF}, \varphi_0)]= \left[({\CF}, \varphi_0) \perp 
\left({\CH}_0({\CQ}_{\bullet})\oplus {\CH}_0({\CQ}_{\bullet})^{\vee}, 
\left(\begin{array}{cc} 0 & 1 \\ 1 &0 \end{array} \right)\right)
\right]
= 
[({\CH}_0({\CU}_{\bullet}), {\CH}_0(\varphi))]=0. 
$$
\pic $\eop$

\section{The Final Results}\label{finalResSec}  
So far we have been working with unsigned translation, in particular in the statement of theorem \ref{divScheme}. 
To conform to the literature,
in this sections we would present our results with respect to the standard signed translation. 
The readers are referred to (\cite{TWGI},
\cite[\S 6]{MSd}) for unexplained notations.
We recall the following notations.
\bE
\item Denote  
$
W_{St}^+({\CA}(X))
:=W\left({\mathcal A}(X), 
^{\vee}, (-1)^{\frac{d(d-1)}{2}}\tilde{\varpi}\right), 
$ and 
\\
$
W_{St}^-({\CA}(X))
:=W\left({\mathcal A}(X), 
^{\vee}, -(-1)^{\frac{d(d-1)}{2}}\tilde{\varpi}\right).  
$
\item Now on, $T:=T_s:D^b({\SV}(X)) \lra D^b({\SV}(X))$ will 
denote the signed ("standard") translation. 
For $n\in {\BZ}$, denote 
$\zeta_n=T^{-n}o\zeta:{\CA}(X)\lra D^b({\SV}(X)$, where $\zeta$ is as in 
(\ref{finalDefZeta}). 
Note $\zeta_0=\zeta$.
\eE 
Now, we state our results.

\bT \label{MainThm}
Suppose $X$ is a noetherian scheme, with $\dim X=d$, as in 
(\ref{notaCate}).
Then,
\bE
\item The functor $\zeta_0:{\mathcal A}(X) \lra 
D^b_{{\mathcal A}(X)}\left({\SV}\left(X\right)\right)$
induces  an  isomorphism
$$
W(\zeta_0): W_{St}^+({\CA}(X)) \iso 
W^d\left(D^b_{{\mathcal A}(X)}({\SV}(X)), *, 1, \varpi \right).
$$
\item The functor $\zeta_1:{\mathcal A}(X) \lra 
D^b_{{\mathcal A}(X)}\left({\SV}\left(X\right)\right)$
induces  an  isomorphism
$$
W_{St}^-({\CA}(X)) \iso W^{d-2}\left(D^b_{{\mathcal A}(X)}
({\SV}(X)), *, 1, -\varpi \right).  
$$
\item For $n=d-1, d-3$, we have
$W^{n}\left(D^b_{{\mathcal A}(X)}({\SV}(X)), *,1, \pm \varpi\right)=0$.  

\eE 
Further, 4-periodicity determines all the shifted Witt groups
$W^{n}\left(D^b_{{\mathcal A}(X)}({\SV}(X)), *,1, \pm \varpi\right)$. 
\eT
\pf Follows from theorem  \ref{divScheme},
as in 
(\cite[\S 6]{MSd}). $\eop$

\noindent{\bf Remark.} 
Following (\cite{MSd}), 
in the case of Cohen-Macaulay affine schemes $X=Spec(A)$, a version of (\ref{MainThm}) was given by Sanders and Sane for resolving subcategories of $Mod(A)$ associated to semi-dualizing modules.

We would also consider the Witt groups of 
the derived category $D^b({\CA}(X))$ and of   
its full subcategory $D^b_{{\CA}}({\CA}(X))$. As was the case with 
$D^b_{{\CA}(X)}({\SV}(X)),
D^b_{{\CA}(X)}({\CA}(X))$ may not have a triangulated structure.
The following follows from the formalism given in (\cite{MSd}).
\bT \label{A2DbA}
Suppose $X$ a noetherian scheme, as in (\ref{notaCate}), with $\dim X=d$.
The duality $^{\vee}:{\CA}\lra {\CA}$ induces a duality
on the derived category $D^b({\CA})$, which we continue to
denote by $^{\vee}$.
Then, $D^b_{{\CA}}({\CA})$ is stable under this duality.
Further, the functors ${\CA}\lra D^b_{{\CA}}({\CA})\hra D^b({\CA})$
induce the following triangle of isomorphisms
$$
\diagram
W({\CA}, ^{\vee}, \pm \tilde{\varpi})\ar[r]^{\sim}\ar[dr]_{\sim} 
& W(D^b_{{\CA}}({\CA}, ^{\vee}, \pm \tilde{\varpi}))\ar[d]^{\wr}\\
&  W(D^b({\CA}, ^{\vee}, \pm \tilde{\varpi}))\\ 
\enddiagram 
$$
\eT
\pf Since ${\CA}$ has 2-out-of-3 property, the
theorem follows from \cite[A.1, A.2]{MSd}.
 $\eop$ 

The following  decomposition of the derived Witt groups is line with the regular
case.
\bT \label{WofProd}
Suppose $X$ is a noetherian scheme, with $\dim X=d$, as in
(\ref{notaCate}) and $X^{(d)}$ will denote the 
set of all closed points of codimension $n$ in $X$.
Then, the  homomorphisms
$$
W^d\left(D^b_{{\mathcal A}}({\SV}(X)), *, 1, \varpi \right)\iso
\bigoplus_{x\in X^{(d)}} 
W^d(D^b_{{\CA}({\CO}_{X,x})}({\CO}_{X,x}),*,1,\varpi) 
\qquad {\rm and} 
$$
$$
W^{d-2}\left(D^b_{{\mathcal A}}({\SV}(X)), *, 1, -\varpi \right)\iso
\bigoplus_{x\in X^{(d)}} 
W^{d-2}(D^b_{{\CA}({\CO}_{X,x})}({\CO}_{X,x}),*,1,-\varpi)
$$
are isomorphims. 
\eT 
\pf We prove only the first isomorphism. 
The diagram:
$$
\diagram 
W_{St}({\CA}(X))^+ \ar[r]^{W(\zeta)\qquad}_{\sim}\ar[d]_{\eta} 
& W^d\left(D^b_{{\mathcal A}(X)}({\SV}(X)), 
*, 1, \varpi \right)\ar[d]^{\gamma}\\
\bigoplus_{x\in X^{(d)}}
W_{St}({\CA}({\CO}_{X,x}))^+ \ar[r]_{\oplus W(\zeta)~~\qquad} 
& \bigoplus_{x\in X^{(d)}} 
W^d(D^b_{{\CA}({\CO}_{X,x})}({\CO}_{X,x}),*,1,\varpi)
\enddiagram 
$$
By theorem (\ref{MainThm}),
two horizontal homomorphisms are isomorphisms.
Also note, 
\linebreak
${\CA}(X) \lra \coprod_{x\in  X^{(d)}} {\CA}({\CO}_{X,x})$
is a duality preserving
equivalence of categories ({\it 
in affine case, this is the Chinese remainder theorem}). 
Therefore, $\eta$ is an isomorphism. Now, since three other maps
in this rectangle are isomorphisms, so is $\gamma$.
\pic $\eop$ 
\begin{remark} \label{regularCoProd}
{\rm
When $X$ is regular, (\ref{WofProd}) is a result of Balmer and 
Walter (\cite{BW}). In this case, 
(\ref{WofProd}) would follow from the
equivalence of the corresponding categories  
$D^b_{\CA}(X)\lra \coprod _{x\in X^{(d)}} {\CA}({\CO}_{X,x})$
(\cite[7.1]{BW}). 
}
\end{remark}
\begin{remark}\label{lowDepth}
{\rm 
It was implicitly assumed in the statement of theorem 
\ref{MainThm} and others, that ${\CA}(X)$ has nonzero objects. 
This would be false, if $depth({\CO}_{X,x})<d=dim X$ for
all closed points $x$. If such cases, a version of these
results can be given by replacing $d$ by $
\delta:=\max\{depth({\CO}_{X,x}): x~{\rm is~a
~closed~point ~of} ~X\}$ (see \S ~\ref{FormalDiv}).
}
\end{remark} 
\section{Witt-D\'{e}vissage formalism}\label{FormalDiv} 
Some authors 
considered resolving subcategories of the category $Mod(A)$ 
of modules over  commutative noetherian rings $A$.
Resolving subcategories of abelian
categories was defined (see \ref{resCatDef}), by substituting the requirement
that $A$ is in resolving subcategory, by the condition 
(3) of the definition 
(\ref{resCatDef}).
In this section, we give  formal versions of
results in \S ~\ref{finalResSec},
 for resolving subcategories
(see \ref{resCatDef}) of abelian categories.
There is a list of examples  of  resolving subcategories of $Mod(A)$ in 
(\cite{T1}). 
The one that is of our particular interest is  
the one corresponding to the notion of semi-dualizing $A$-modules $\omega$.
Such semi-dualizing $A$-modules $\omega$ 
naturally give rise to  resolving
subcategories  of $Mod(A)$.
With that in mind, we define
the following and  establish our set up for this section.

\bD\label{defResWduality}{\rm 
Suppose ${\SV}$ is a resolving
subcategory of an abelian category ${\CC}$. 
Let $\omega$ be a fixed object on ${\CC}$ and write 
$M^*:=Mor(M, \omega)$. Assume that $\omega$ has an injective
resolution, and choose one such resolution ${\SI}_{\bullet}$,  as follows:
\begin{equation}  
\diagram
0\ar[r] & \omega \ar[r] & {\mathscr I}_0 \ar[r] & {\SI}_1\ar[r] 
& {\SI}_2\ar[r] & \cdots
\enddiagram 
\end{equation}  
\bE
\item 
Given any object $M$ in ${\CC}$, define $Ext^i(M, \omega):=H^i(M, 
{\SI}_{\bullet})$. 
({\it We only define $Ext^i(-, \omega)$; and
not giving any definition
of $Ext^i(-, {\CG})$ when ${\CG}\neq \omega$.}) It immediately follows that 
$\{M\mapsto Ext^i(-, \omega): i\geq 0\}$ is a contravariant $\delta$-functor
(see \cite[\S III.1]{H}). 

\item 
We say that ${\SV}$ inherits 
a {\bf $\omega$-duality structure} if, 
for all objects $P$ in ${\SV}$,
(1) $P^*\in {\SV}$, 
(2) and there is a natural equivalence $ev: P\iso P^{**}$,
from the identity fucntor $1_{{\SV}}$ to the double dual.
\item Also, we say that ${\SV}$ is  {\bf totally  $\omega$-reflexive} if
in addition, 
for all objects $P$ in ${\SV}$ and integers $i\geq 1$, 
$Ext^i(P, \omega)=0$.
It is easy to see that if $P_{\bullet}\sur M$ is a ${\SV}$-resolution 
of an object in ${\CC}$ that there are natural isomorphisms 
from the homologies 
$H_i(P_{\bullet}^{\#})\iso Ext^i(M, \omega)$, where $P_{\bullet}^{\#}$
denote the dual of $P_{\bullet}$ induced by $^*$. 

\item We also recall some standard definitions.
For objects $M$ in ${\CC}$, $\dim_{{\SV}}M$ will denote the 
length of the shortest ${\SV}$-resolution of $M$ (which can be infinite).
 As in \S ~\ref{resZetaSec}, 
${\CB}:={\CB}({\SV}):=\{M\in {\CC}: \dim_{{\SV}}M<\infty \}$ denotes 
the full subcategory of such objects.
 Also denote, 
$d:=\dim_{{\SV}}({\CB}):=\max\{ \dim_{{\SV}}M: M\in {\CB} \}$. 
If $d=\dim_{{\SV}}({\CB})<\infty$, let 
$$
{\CA}(\omega):= 
{\CA}({\SV}, \omega):=\{M \in {\CB}({\SV}): Ext^i(M,\omega)=0 ~~\forall ~i<d. 
\}
$$
denote the full subcategory. 
It follows ${\CB}({\SV})$ and  
${\CA}({\omega})$ are exact subcategories of ${\CC}$. 
\item \label{setUp4Res} {\bf Set up}:
In what follows, we will have the following set up:

\noindent{\it ${\CC}$ will denote an abelian category and ${\SV}$ will be a
resolving subcategory of ${\CC}$. We fix an object $\omega$ with 
a (chosen) injective resolution $\omega \hra {\SI}_{\bullet}$.
The $\delta$-functor $\{Ext^i(-,\omega)\}$ is defined as above.  
We assume:
(1) ${\SV}$ inherits a $\omega$-duality structure and is totally
$\omega$-reflexive.
(2) Given an object ${\CF}\in {\CB}$ and a
${\SV}$-resolution ${\CE}_{\bullet}\sur {\CF}$,
the cycle objects
$Z_n:=\ker({\CE}_{n}\lra {\CE}_{n-1})\in {\SV}$~~$\forall~n\gg 0$.
(3) Further, $d:=\dim_{{\SV}}({\CB})<\infty$. 
}
We observed (\ref{formalZeta}), under this set up 
\bE
\item 
 ${\CB}(\SV)$ is an exact subcategory and every epimorphism in
${\CB}$ is admissible.
\item \label{2OfForZeta}
Further, the  functor
$\zeta: {\CB}({\SV}) \lra D^b(\SV)$
is defined (see \ref{defFunt(g)}).
\item It also follows that ${\CA}(\omega)$ is an exact subcategory of 
${\CC}$. 
\eE

\eE 
}\eD 
Most of what are in  \S ~\ref{finalResSec}, work
for resolving subcategories ${\SV}$, of  abelian categories ${\CC}$,
as in the set up (\ref{setUp4Res}) of \ref{defResWduality}.  
We state them below.
\bL\label{dualRESCheck}
Suppose 
$({\CC}, {\SV}, \omega)$ is as in (\ref{setUp4Res}) of  \ref{defResWduality}.
Then, the  association $M\mapsto {M}^{\vee}:= Ext^d\left(M, \omega \right)$
defines 
a duality $^{\vee}:{\CA}(\omega) \lra {\CA}(\omega)$. 
 \eL
\pf Note,  $M\in {\CA}(\omega)\Lra  M^{\vee}\in {\CA}(\omega)$.  
The rest of the proof is  as that of (\ref{dualCheck}). $\eop$

\bT\label{dualResOfDAV}
Suppose
$({\CC}, {\SV}, \omega)$ be as in (\ref{setUp4Res}) of  \ref{defResWduality}. Suppose ${\CE}_{\bullet}$ is a complex in $Ch^b_{{\CA}}({\SV})$.
Then, the dual ${\CE}_{\bullet}^{\#}$ is also $Ch^b_{{\CA}}({\SV})$.
Further, there is a canonical isomorphism $\eta: H_{-r}({\CE}_{\bullet}^{\#})
\iso H_{r-d}({\CE}_{\bullet})^{\vee}$, which is natural with respect to 
morphisms  in $Ch^b_{{\CA}}({\SV})$.   
\eT 
\pf Same as that of (\ref{dual4ChbPX}). $\eop$

\bT \label{MainRESThm}
Suppose $({\CC}, {\SV}, \omega)$ be as in (\ref{setUp4Res}) of  \ref{defResWduality}
and the rest 
of the notations  be same as in \S ~\ref{finalResSec}. 
Also refer  to the definition (\ref{formalZeta})
of the functor $\zeta: {\CB}({\SV})
\lra D^b({\SV})$.  
Then,
\bE
\item The functor $\zeta_0:{\mathcal A}(\omega) \lra 
D^b_{{\mathcal A}(\omega)}\left({\SV}\right)$
induces  an  isomorphism
$$
W(\zeta_0): W_{St}^+({\CA}(\omega)) \iso 
W^d\left(D^b_{{\mathcal A}(\omega)}({\SV}), *, 1, \varpi \right).
$$
\item The functor $\zeta_1:{\mathcal A}(\omega) \lra 
D^b_{{\mathcal A}(\omega)}\left({\SV}\right)$
induces  an  isomorphism
$$
W_{St}^-({\CA}(\omega)) \iso W^{d-2}\left(D^b_{{\mathcal A}(\omega)}
({\SV}), *, 1, -\varpi \right).  
$$
\item For $n=d-1, d-3$, we have
$W^{n}\left(D^b_{{\mathcal A}(\omega)}({\SV}), *,1, \pm \varpi\right)=0$.

\eE
Further, 4-periodicity determines all the shifted Witt groups
$W^{n}\left(D^b_{{\mathcal A}(\omega)}({\SV}), *,1, \pm \varpi)\right)$.
\eT
\pf 
Similar to that of theorem \ref{MainThm}. 
$\eop$

\bT \label{A2DbARES}
Suppose
$({\CC}, {\SV}, \omega)$ is as in (\ref{setUp4Res}) of  \ref{defResWduality}.
Then, $D^b_{{\CA}(\omega)}({\CA}(\omega))$ is stable under the duality
on $D^b({\CA}(\omega))$
induced by $^{\vee}:{\CA}(\omega)\lra {\CA}(\omega)$
{\rm (which we continue to
denote by $^{\vee}$)}.
Further, the functors
${\CA}(\omega)\lra D^b_{{\CA}(\omega)}({\CA}(\omega))\hra D^b({\CA}(\omega))$
induce the following triangle of isomorphisms
$$
\diagram
W({\CA}(\omega), ^{\vee}, \pm \tilde{\varpi})\ar[r]^{\sim}\ar[dr]_{\sim} 
& W(D^b_{{\CA}(\omega)}({\CA}(\omega), ^{\vee}, \pm \tilde{\varpi}))\ar[d]^{\wr}\\
&  W(D^b({\CA}(\omega), ^{\vee}, \pm \tilde{\varpi}))\\ 
\enddiagram 
$$
\eT
\pf Note that the diagonal isomorphism follows
directly from (\cite[4.7]{TWGII}). So, we need only to prove that
the horizontal map is surjective.
The proof would be very similar to that of (\ref{A2DbA}),
by an application \cite[A.1, A.2]{MSd}.
Some clarification is needed, because ${\CA}(\omega)$
does not have the 2-out-of-3 property, which was used in
(\ref{A2DbA}).
%
First, we claim that for complexes $({\CE}_{\bullet},
\partial_{\bullet})$ in
$Ch^b_{{\CA}(\omega)}({\CA}(\omega))$ all the boundaries
$B_i:=image(\partial_{i-1})\subseteq {\CE}_i$
and  cycles $Z_i:=\ker(\partial_i)\subseteq {\CE}_i$ are in ${\CA}(\omega)$.
To see this note, since ${\CB}({\SV})$ is exact and epimorphisms
in ${\CB}({\SV})$ are admissible, inducting from right, it follows
$B_i, Z_i, \frac{{\CE}_i}{B_i}$ are in ${\CB}({\SV})$.
Now, for
$X= B_i, Z_i, \frac{{\CE}_i}{B_i}$, by induction from left, it follows
$Ext^i(X, \omega)=0$
for all $i\neq d$, hence
are in ${\CA}(\omega)$.
This establishes the claim. In fact, proof of  \cite[A.1, A.2]{MSd}
works whenever, 
$B_i, Z_i, \frac{{\CE}_i}{B_i}$
are in ${\CA}(\omega)$ for all $i$.
This completes the proof.

We comment that one  can combine \ref{formalZeta} and \ref{dualResOfDAV}
to give a proof of the duality statement,
by constructing a double complex in $D^b({\SV})$,
as in
(\cite[\S 3]{MSd}).
$\eop$

\end{document}